\documentclass{scrartcl}
\usepackage[numbers]{natbib}
\usepackage{fontenc}
\usepackage{shadethm}
\usepackage{amsthm}
\usepackage{float}
\usepackage{bm}
\usepackage{graphicx}
\usepackage{subfig}
\usepackage{stmaryrd}
\usepackage{mathrsfs} 
\usepackage{amsmath}
\usepackage{amssymb}
\usepackage{wasysym}
\usepackage{sectsty}
\usepackage{hyperref}
\usepackage[a4paper, left=2.7cm, right=2.7cm, top=2.5cm]{geometry}
\usepackage{chngcntr}
\counterwithin*{equation}{section}

\usepackage{cleveref}
\DeclareMathAlphabet{\mathpzc}{OT1}{pzc}{m}{it}

\sectionfont{\normalsize\centering}
\subsectionfont{\footnotesize\centering}

\theoremstyle{plain}
\newtheorem{thm}{Theorem}[section] 

\theoremstyle{definition}
\newtheorem{defn}[thm]{Definition} 
\newtheorem{lem}[thm]{Lemma}
\newtheorem{prop}[thm]{Proposition}
\newtheorem{rem}[thm]{Remark}

\def\XXint#1#2#3{{\setbox0=\hbox{$#1{#2#3}{\int}$ }
		\vcenter{\hbox{$#2#3$ }}\kern-.6\wd0}}

\usepackage[utf8]{inputenc}
\usepackage[german,english,russian]{babel}

\newcounter{MPequ}

\newcounter{AppA}

\begin{document}\selectlanguage{english}
\begin{center}
\normalsize \textbf{\textsf{Isoperimetric problem for the first curl eigenvalue}}
\end{center}
\begin{center}
	Wadim Gerner\footnote{\textit{E-mail address:} \href{mailto:wadim.gerner@icmat.es}{wadim.gerner@icmat.es}}
\end{center}
\begin{center}
{\footnotesize	Instituto de Ciencias Matem\'{a}ticas, Consejo Superior de Investigaciones Cient\'{i}ficas, 28049 Madrid, Spain}
\end{center}
{\small \textbf{Abstract:} We consider an isoperimetric problem involving the smallest positive and largest negative curl eigenvalues on abstract ambient manifolds, with a focus on the standard model spaces. We in particular show that the corresponding eigenvalues on optimal domains, assuming optimal domains exist, must be simple in the Euclidean and hyperbolic setting. This generalises a recent result by Enciso and Peralta-Salas who showed the simplicity for axisymmetric optimal domains with connected boundary in Euclidean space.
\newline	
We then generalise another recent result by Enciso and Peralta-Salas, namely that the points of any rotationally symmetric optimal domain with connected boundary in Euclidean space which are closest to the symmetry axis must disconnect the boundary, to the hyperbolic setting, as well as strengthen it in the Euclidean case by getting rid of the connected boundary assumption.
\newline
Lastly, we show how a second variation inequality related to the isoperimetric problem may be used in order to relate the existence of Killing-Beltrami fields to the geometry of the ambient space.
\newline
\newline
{\small \textit{Keywords}: Isoperimetric problems, Spectral theory, Beltrami fields, Killing fields, Curl operator}
\newline
{\small \textit{2020 MSC}: 35P15, 35Q31, 35Q35, 35Q85, 58C40, 76W05}
\section{Introduction}
Eigenvector fields of the curl operator, so called (strong) Beltrami fields appear naturally in physics, most notably in magnetohydrodynamics and fluid dynamics. On the one hand they appear as special solutions of the incompressible Euler equations \cite{A66} and on the other hand as stationary solutions of the equations of ideal magnetohydrodynamics in the context of a resting plasma and constant pressure \cite[Chapter III]{AK98}. In fluid dynamics, Beltrami fields are of particular interest from a topological perspective. While their non-Beltrami relatives in essence behave all the same \cite{A66}, a fact known as Arnold's structure theorem, the behaviour of Beltrami flows can be much more delicate, see for instance \cite{DFGHMS86}, \cite{EP12}, \cite{EP-S15} and \cite{EPT17}. This begs the question of whether or not on any given (compact) manifold we can always find Beltrami fields. It was observed by Woltjer \cite{W58} in the context of ideal magnetohydrodynamics that there is a variational formulation to this problem. Namely, Woltjer considered formally a Lagrange-multiplier approach in the presence of the so called helicity constraint. Helicity turns out to be formally preserved by the evolution equations of ideal magnetohydrodynamics and a physical interpretation of helicity was provided by Moffatt \cite{M69}. It was then first rigorously shown by Arnold by means of a spectral theoretical argument, \cite{A74b}, that this variational problem admits solutions in the setting of closed manifolds and that the solutions are, as expected, Beltrami fields. Solutions to this variational problem in the context of bounded domains with smooth boundary in Euclidean space were then obtained for example in \cite{AL91} and \cite{CDG00} and a detailed discussion in the context of abstract compact manifolds with boundary can be found in \cite{G20}. It follows from \cite[Theorem 2.1]{G20}, but is also implicitly contained in the works \cite{A74b} and \cite{CDG00} for their respective setting, that the solutions to the variational constraint minimisation problem posed by Woltjer are precisely those vector fields which are eigenfields of the curl operator corresponding to its smallest positive or largest negative eigenvalue, see also \cref{MT1} of the present paper for a precise statement regarding the existence of these eigenvalues. Further, one can show that the helicity functional, which appears as the constraint in the variational formulation, is invariant under the action of volume preserving diffeomorphisms, see \cite[Section 2.3 Corollary]{A74b} and \cite[Lemma 4.5]{G20}. In fact it turns out that helicity is in essence the only such invariant \cite{EPT16}, see also \cite{KPY20}. In view of this it is then natural to ask whether or not it is possible to minimise the first curl eigenvalue in fixed volume classes for some given reference ambient manifold, such as $\mathbb{R}^3$. Let us for simplicity focus for now on the $\mathbb{R}^3$ case, even though the results we obtain are more general. We observe that our curl eigenfields are in particular divergence-free so that it follows from standard calculus identities that any such eigenfield $\bm{X}$ also satisifes
\[
-\Delta \bm{X}=\text{curl}(\text{curl}(\bm{X}))=\lambda_+^2\bm{X},
\]
with $\lambda_+$ being our smallest positive curl eigenvalue and $\Delta$ denoting the standard Euclidean Laplacian. Hence, the smallest positive curl eigenvalue turns out to be a Laplace-eigenvalue. The isoperimetric problem for the first Dirichlet-Laplace-eigenvalue has a long history and it was conjectured by Lord Rayleigh that the first eigenvalue is always minimised by the round ball of the same volume, \cite{LR94} (see \cite{LR45} for a reprinted version). This conjecture has been resolved by means of the Faber-Krahn inequality, see the original references \cite{Fa23} and \cite{Kr25} which address the $\mathbb{R}^2$ case. While there exist some Faber-Krahn type inequalities for the first curl eigenvalue problem in different contexts, see \cite[Theorem B]{CDG00}, \cite[Theorem 5.6]{P04} and recent work of Enciso and Peralta-Salas \cite[Appendix A, Theorem A.1]{EPS21}\footnote[2]{arXiv identifier: 2007.05406}, in the sense that the first curl eigenvalue may be bounded away from zero by some function, which only depends on the volume of the domain, the question of existence of optimal domains in the context of the first curl eigenvalue is still wide open. Interestingly enough not many results may be found in the literature regarding this curl isoperimetric problem. A helicity maximisation problem, which is somewhat related to (but distinct from) the isoperimetric problem under consideration here, was studied in \cite{CDGT002} in the Euclidean setting. The curl isoperimetric problem on general ambient spaces was studied by the present author in his recent PhD thesis where also a second variation inequality was derived \cite[Chapter 2]{G20Diss}. It was conjectured in \cite{CDGT002} that their corresponding helicity isoperimetric problem, which is related to, but distinct from, the problem addressed in the present work, does not admit a solution. A partial result in this direction, in the context of the curl isoperimetric problem studied in the present paper, was obtained in \cite{EPS21} by Enciso and Peralta-Salas concerning the non-optimality of a certain class of rotationally symmetric domains. Interestingly enough, the arguments they provide do not generalise to the setting of Cantarella et al. in \cite{CDGT002}. The main reason being that in the setting of \cite{CDGT002} the considered vector fields are in general not $L^2$-orthogonal to the space of harmonic vector fields.
\newline
One main step in the work of Enciso and Peralta-Salas is that they show the simplicity of the first curl eigenvalue on rotationally symmetric optimal domains with connected boundary (assuming such domains exist). The first result of the present work, \cref{MT2}, generalises this statement to the setting of arbitrary optimal domains. The arguments used in the present work, which allow for this generalisation, lie in the observation that $L^2$-orthogonality of curl eigenfields translates on optimal domains to pointwise $g$-orthogonality on the boundary, with $g$ being our Riemannian metric. Note that contrary to the situation of the first Dirichlet-Laplacian eigenvalue, the first curl eigenvalue can in general have multiplicity strictly larger than $1$, \cite{CDGT00}, so that these upper bounds are far from trivial.
\newline
In \cref{MC8} we generalise the result obtained by Enciso and Peralta-Salas in Euclidean space \cite{EPS21} to the setting of rotationally symmetric domains with not necessarily connected boundary. The main ingredient here is the newly obtained simplicity of the first curl eigenvalue for all optimal domains and the fact that we can associate an appropriate, rotationally symmetric vector potential to each corresponding first curl eigenfield.
\newline
In the last part of this article we establish a connection between the existence of so called Killing-Beltrami fields, i.e. eigenfields of curl whose flows induce isometries, and the geometry of the ambient space, namely its sectional curvature. The ideas of the proof are closely related to the derivation of the second variation inequality in the context of the isoperimetric problem considered throughout the present work.
\newline
Let us finally stretch one more time that the helicity maximisation problem studied in \cite{CDGT002} is distinct from the spectral problem studied in the present work and that there is no reason to believe that the optimal domains in these two optimisation problems coincide. The study of the spectral problem, as considered in the present work, was initiated by Enciso and Peralta-Salas\footnote[3]{Preprint published on 10th July 2020} in \cite{EPS21} and independently by the present author in his dissertation\footnote[4]{Successfully defended on 10th December 2020} \cite[Chapter 2]{G20Diss}.
\section{Main results}
\textbf{Conventions:} We will always consider an \textit{ambient space} $(\mathcal{R},g)$, which is assumed to be a connected, oriented, smooth Riemannian $3$-manifold without boundary. We will simply say "Let $(\mathcal{R},g)$ be our ambient space" by which we assume the previously mentioned properties. We refer to smoothly embedded, connected, compact $3$-manifolds $\bar{M}\subseteq \mathcal{R}$ with boundary as \textit{compact domains} of the ambient space and they are always equipped with the induced structures of the ambient space. We let $\text{Sub}_c(\mathcal{R})$ denote the set of all compact domains $\bar{M}$ of some given ambient space $(\mathcal{R},g)$ and for given $0<V<\text{vol}(\mathcal{R})$ we let $\text{Sub}^V_c(\mathcal{R})\subset \text{Sub}_c(\mathcal{R})$ be those compact domains $\bar{M}$ with $\text{vol}(\bar{M})=V$.
\newline
We denote by $\mathcal{V}(\bar{M})$ the set of all smooth vector fields on a manifold $\bar{M}$, by $\mathcal{V}_n(\bar{M})$ the set of all smooth vector fields $\bm{X}$ such that there exists some $\bm{A}\in \mathcal{V}(\bar{M})$ with $\bm{A}\perp \partial\bar{M}$ and $\text{curl}(\bm{A})=\bm{X}$.
\newline
We refer to a vector field $\bm{X}\in \mathcal{V}_n(\bar{M})$ as a (strong) Beltrami field if it is an eigenvector field of the curl operator corresponding to a non-zero eigenvalue (in particular $\bm{X}\not\equiv 0$). Note that the notion of Beltrami fields is usually more general in the literature, see for example \cite{N14}.
\newline
We use the common notions of divergence and curl of vector fields on abstract manifolds, i.e. $\text{div}(\bm{X}):=-\delta \omega^1_{\bm{X}}$, where $\delta$ denotes the adjoint derivative and $\omega^1_{\bm{X}}$ the associated one form with $\bm{X}$ by means of some fixed Riemannian metric $g$, and $\omega^1_{\text{curl}(\bm{X})}=\star d\omega^1_{\bm{X}}$, where $\star$ denotes the Hodge-star operator and $d$ denotes the exterior derivative.
\newline
Lastly, we denote by $\mathbb{R}^3$ the standard Euclidean $3$-space, by $\mathcal{H}^3$ the hyperbolic $3$-space of constant curvature $-1$, by $S^3$ the round $3$-sphere of constant curvature $+1$ and by $S^3_+:=\{(x,y,z,w)\in S^3\subset \mathbb{R}^4|w>0\}$ the upper half sphere with the same metric.
\newline
\newline
Let us start by quoting the following result
\begin{thm}[\cite{A74b}, \cite{G20}]
	\label{MT1}
Let $(\bar{M},g)$ be a compact, oriented, smooth Riemannian $3$-manifold with (possibly empty) boundary. Then the curl operator
\[
\text{curl}:\mathcal{V}_n(\bar{M})\rightarrow \mathcal{V}(\bar{M}),\text{ }\bm{X}\mapsto \text{curl}(\bm{X})
\]
admits a smallest positive eigenvalue $\lambda_+(\bar{M})>0$ and a largest negative eigenvalue $\lambda_-(\bar{M})<0$.
\end{thm}
Given some ambient space $(\mathcal{R},g)$ and some $0<V<\text{vol}(\mathcal{R})$, we are interested in the following minimisation problem
\begin{align}
	\label{ME1}
	\text{Sub}^V_c(\mathcal{R})\rightarrow \mathbb{R},\text{ }\bar{M}\mapsto \lambda_+(\bar{M})\rightarrowtail\text{ Min. }
\end{align}
The existence of minimisers of the above problem is still an open problem. In the following we will refer to minimisers as \textit{optimal domains}. Here we focus on the eigenvalue $\lambda_+(\bar{M})$, but the same proofs apply equally well to the corresponding setting of $\lambda_-(\bar{M})$.
\newline
Regarding our standard model spaces $\mathbb{R}^3$, $\mathcal{H}^3$ and the half sphere $S^3_+$ we have the following main result concerning the multiplicity of the first eigenvalue, where we set $E_+(\bar{M}):=\{\bm{X}\in \mathcal{V}_n(\bar{M})|\text{curl}(\bm{X})=\lambda_+(\bar{M})\bm{X} \}$.
\begin{thm}[Multiplicity of the eigenvalue]
	\label{MT2}
	Let $(\mathcal{R},g)\in \{\mathbb{R}^3,\mathcal{H}^3,S^3_+ \}$ be our ambient space. Given some $0<V<\text{vol}(\mathcal{R})$, if there exists an optimal domain $\bar{M}\in \text{Sub}^V_c(\mathcal{R})$, then $\lambda_+(\bar{M})$ is simple, i.e. $\dim\left(E_+(\bar{M})\right)=1$.
\end{thm}
Note that on general compact domains the eigenvalue need not be simple \cite{CDGT00}. Further, the simplicity of the eigenvalue in the connected boundary case was established for the Euclidean $3$-space in \cite{EPS21} under the additional assumption of axisymmetry. Their proof, however, does not extend to general non-axisymmetric, compact domains. Hence, \cref{MT2} may be viewed as an extension of this result.
\newline
\newline
The situation for the remaining model space $S^3$ turns out to be more complicated.
\begin{thm}[Eigenvalue $\lambda_+$ on $S^3$]
	\label{MT3}
	Let $S^3$ be our ambient space and let $0<V<\text{vol}(S^3)=2\pi^2$ be given. If there exists an optimal domain $\bar{M}\in \text{Sub}^V_c(S^3)$ with connected boundary such that its first eigenvalue $\lambda_+(\bar{M})$ is not simple, then (after possibly applying isometries to $\bar{M}$) the function
	\[
	f:S^3\subset \mathbb{R}^4\rightarrow \mathbb{R},\text{ }(x,y,z,w)\mapsto w
	\]
	must integrate to $0$ over $\bar{M}$, i.e. $\int_{\bar{M}}f\omega_g=0$, where $\omega_g$ denotes the "round" volume form on $S^3$.
\end{thm}
\Cref{MT2}, \cref{MT3} and the discrepancy between the model spaces are rooted in the existence of so called concircular vector fields, see \cref{S2D5}, and the behaviour of the corresponding potential functions.
\newline
\newline
The upcoming result deals with the question of the existence of rotationally symmetric domains in hyperbolic $3$-space. To this end we define the following vector field
\begin{align}
	\label{ME4}
	\bm{R}:\mathbb{R}^3\rightarrow \mathbb{R}^3,\text{ }(x,y,z)\mapsto (-y,x,0),
\end{align}
which induces rotations around the $z$-axis and we consider the Poincar\'{e} ball model of hyperbolic $3$-space. Then $\bm{R}$ gives rise to a flow of isometries. The natural question in this context is whether or not a rotationally symmetric optimal domain can exist, i.e. whether there exists some optimal $\bar{M}\in \text{Sub}^V_c(\mathcal{H}^3)$ for some $0<V<+\infty$ such that $\bm{R}\parallel \partial\bar{M}$. This problem was studied in \cite{EPS21} in the context of Euclidean space. The following is the corresponding result in hyperbolic space.
\begin{thm}[Rotationally symmetric optimal domains]
	\label{MT5}
	Let $\mathcal{H}^3$ be our ambient space and $0<V<+\infty$ be given. Let further $g_{\mathcal{H}}$ denote the corresponding hyperbolic metric in the Poincar\'{e} ball model. If there exists an optimal domain $\bar{M}\in \text{Sub}^V_c(\mathcal{H}^3)$ such that $\bm{R}\parallel \partial\bar{M}$, see (\ref{ME4}), then $\bar{M}$ does not intersect the $z$-axis. Further, if $\partial\bar{M}$ is connected, the set $N_0:=\left\{p\in \partial\bar{M}| |\bm{R}(p)|_{g_{\mathcal{H}^3}}=\min_{q\in \partial\bar{M}}|\bm{R}(q)|_{g_{\mathcal{H}^3}} \right\}$ must disconnect the boundary, i.e. the non-empty set $\partial \bar{M}\setminus N_0$ must be disconnected.
\end{thm}
Here $|\bm{R}(q)|_{g_{\mathcal{H}^3}}=\sqrt{g_{\mathcal{H}^3}(\bm{R}(q),\bm{R}(q))}$. As a consequence, just like in the Euclidean case, no rotationally symmetric, compact domain $\bar{M}\in \text{Sub}_c(\mathcal{H}^3)$ with connected boundary and for which $\partial\bar{M}\setminus N_0$ is connected can be optimal.
\begin{rem}
	\label{MR6}
	By means of the stereographic projection we can identify $(S^3_+,g_R)$ equipped with the round metric $g_R$ with the open Euclidean unit ball $B_1(0)\subset \mathbb{R}^3$ equipped with the metric $g(\vec{x})=\frac{4}{(1+|\vec{x}|_{2}^2)^2}g_E(\vec{x})$, where $\vec{x}=(x,y,z)\in B_1(0)$, $|\cdot|_{2}$ denotes the Euclidean distance and $g_E$ denotes the standard Euclidean metric. Then the vector field $\bm{R}$ gives also rise to a Killing field on $S^3_+$, so that one can consider the corresponding problem of \cref{MT5} on $S^3_+$ as well.
\end{rem}
In the setting of $\bm{R}\parallel \partial\bar{M}$ one can in fact say even more. The upcoming result is an improvement of the main result in \cite{EPS21}, since we do not demand connectedness of the boundary. Here we focus on the standard Euclidean situation.
\newline
Before we state the result, we need to introduce some notation. We will support the given mathematical definitions with suitable figures, see \cref{Fig1}. Just like in \cref{MT5} we will see that any rotationally symmetric domain in Euclidean space, which intersects the $z$-axis, is not optimal. Therefore we may consider a cross section $C\cong (S^1)^{\#\partial\bar{M}}$, where $\#\partial\bar{M}$ denotes the number of connected components of $\partial\bar{M}$, in the $x$-$z$-plane, where each boundary component is represented by some closed smoothly embedded curve and such that $\partial\bar{M}\cong S^1\times C$. There is a unique closed curve $\gamma_0$ whose image contains the points of the boundary in the cross section which are closest to the $z$-axis. In the following we will denote by $\gamma_0$ the curve, as well as its image. In correspondence with \cref{MT5} we define
\begin{align}
	\label{ME7}
	d_-:=\min_{q\in \gamma_0}\{|\bm{R}(q)|_{2}\}=\min_{q\in \partial\bar{M}}\{|\bm{R}(q)|_{2}\},\text{ } N_0:=\{p\in \gamma_0| |\bm{R}(p)|_{2}=d_-\},
\end{align}
where we recall that $|\cdot|_2$ denotes the standard Euclidean metric. It then follows by compactness that there are points $\vec{x}_+,\vec{x}_-\in N_0$ for which the $z$-component becomes maximal, respectively minimal. Then the points $\vec{x}_+$ and $\vec{x}_-$ divide the curve $\gamma_0$ in two parts. Precisely one of these parts will contain those points of $\gamma_0$ which have maximal distance from the $z$-axis and a second curve which does not contain such points. We denote by $L_-$ the part of $\gamma_0$ connecting $\vec{x}_+$ with $\vec{x}_-$ and \underline{not} containing points of maximal distance to the $z$-axis. Further, if the boundary has more than one connected component, there will be smoothly embedded, closed curves $\gamma_1,\dots,\gamma_N$ with $N=\#\partial\bar{M}-1$, representing the remaining boundary components. Again by compactness, there must be a point in $L_-\cup \cup_{i=1}^N\gamma_i$ which maximises the distance to the $z$-axis on this set. We denote the maximal distance on this set by $d_+$. We finally let $L_+\subseteq \gamma_0\setminus L_-$ denote the part of $\gamma_0$ which has a distance greater or equal to $d_+$ to the $z$-axis, i.e.
\begin{align}
	\label{ME8}
	L_+:=\{p\in \gamma_0\setminus L_-| |\bm{R}(p)|_{2}\geq d_+ \}.
\end{align}
We note that $L_+$ is always non-empty. For an illustration of the situation see \cref{Fig1}.
\begin{thm}[Rotationally symmetric optimal domains in $\mathbb{R}^3$]
	\label{MC8}
	Let $\mathbb{R}^3$ be our ambient space and $\bar{M}\in \text{Sub}_c(\mathbb{R}^3)$, such that $\bm{R}\parallel \partial\bar{M}$. If either $\bar{M}$ intersects the $z$-axis or otherwise if $\mathcal{L}^1(L_+)\geq \mathcal{L}^1(L_-)+\sum_{i=1}^N\mathcal{L}^1(\gamma_i)$, then $\bar{M}$ is not an optimal domain in its own volume class, where $N=\#\partial\bar{M}-1$.
\end{thm}
Here $\mathcal{L}^1$ denotes the $1$-dimensional Hausdorff measure, which coincides with the length of the respective curves in the sense of Riemannian geometry. Further, if $\partial\bar{M}$ is connected, then the sum $\sum_{i=1}^N\mathcal{L}^1(\gamma_i)$ is set to zero by convention.
\begin{figure}
	\centering
	\subfloat[Here $N_0$ consists of exactly the two green points. The closed purple curve describes a second boundary component.]{
		\centering
		\includegraphics[width=0.35\textwidth, keepaspectratio]{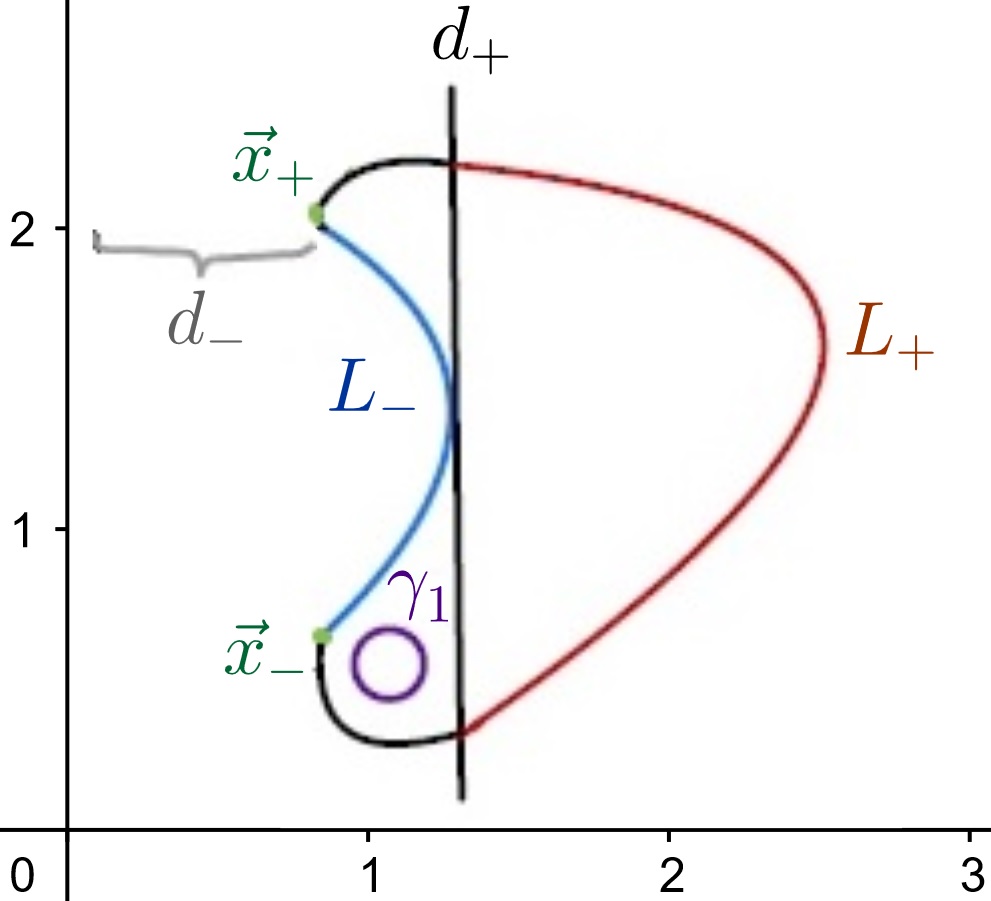}
		\label{Fig1a}
	}
	\hspace{1cm}
	\subfloat[Here $N_0$ consists of $3$ points. $\gamma_0$ is the closed outer smooth curve consisting of the red, blue and black arcs.]{
		\centering
		\includegraphics[width=0.5\textwidth, keepaspectratio]{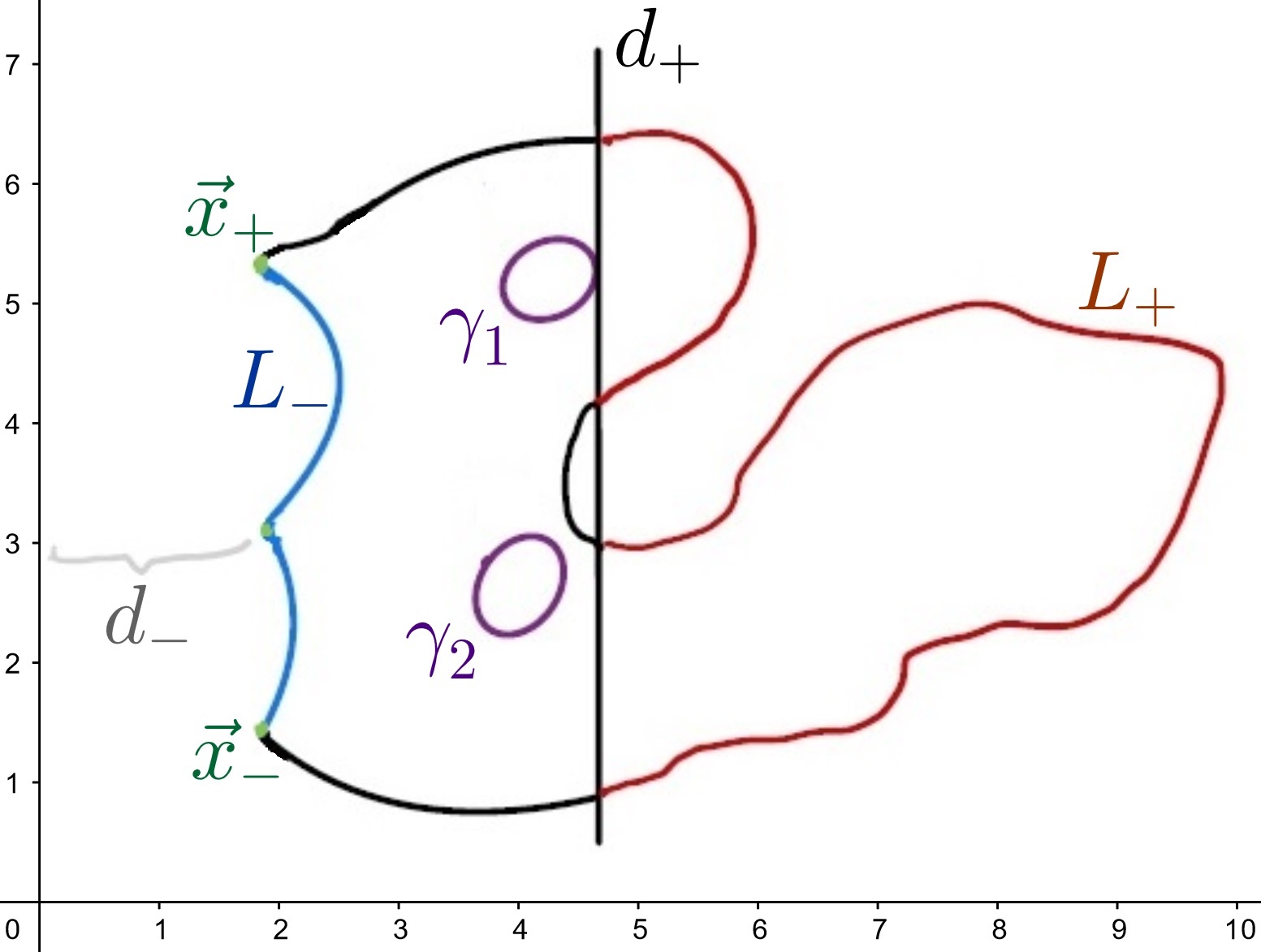}
		\label{Fig1b}
	}
	\caption{Illustration of \cref{MC8}.}
	\label{Fig1}
\end{figure}
\newline
For an example of a cross section $C$, which does not give rise to an optimal domain by means of \cref{MC8}, see \cref{Fig1a} and \cref{Fig1b}. In both cases the length of the red curves is greater or equal to the sum of the lengths of the blue curve and the purple curves. Hence a domain obtained from the given cross sections by rotation around the $z$-axis won't give rise to optimal domains.
\newline
Lastly note that the analogue of \cref{MT5} for the Euclidean setting follows from \cref{MC8}, since if $\partial\bar{M}$ and $\partial\bar{M}\setminus \{p\in \partial\bar{M}||\bm{R}(p)|_{g_E}=d_- \}$ are both connected, it follows that $N_0$ must be connected and hence must be diffeomorphic to a closed interval (possibly a point). By its definition $L_-$ must then be the line segment parallel to the $z$-axis connecting $\vec{x}_+$ and $\vec{x}_-$. On the other hand we find in this scenario $L_+=\gamma_0\setminus L_-$ and we know that this curve also connects the points $\vec{x}_+$ and $\vec{x}_-$. But since $L_-$ is a distance minimising geodesic, we must have $\mathcal{L}^1(L_+)\geq \mathcal{L}^1(L_-)$ and therefore no such domain can be optimal.
\newline
\newline
\Cref{MT2}-\cref{MC8} are the main results of the present work. In the appendix we included a discussion on how a second variation approach connected to our optimal domain problem allows us to establish a connection between the existence of Killing-Beltrami fields and the geometry of the ambient space.
\section{Proof of \cref{MT2} and \cref{MT3}}
\subsection{Concircular vector fields and divergence-free extensions}
Let us first start with a definition
\begin{defn}[div-free extension]
	\label{S2D1}
	Let $(\mathcal{R},g)$ be our ambient space. We say that $(\mathcal{R},g)$ has the div-free extension property, if for every $\bar{M}\in \text{Sub}_c(\mathcal{R})$ and every $\bm{X}\in \mathcal{V}(\bar{M})$ with $\text{div}(\bm{X})=0$, there exists some open subset $\bar{M}\subseteq U\subseteq \mathcal{R}$ and a divergence-free extension $\widetilde{\bm{X}}\in \mathcal{V}(U)$ of $\bm{X}$ to $U$.
\end{defn}
For our purposes the following results will be enough, where we define for given $\bar{M}\in \text{Sub}_c(\mathcal{R})$
\begin{gather}
	\label{S2E1}
	\mathcal{H}_{\text{co}}(\bar{M}):=\{\bm{\Gamma}\in \mathcal{V}(\bar{M})|\bm{\Gamma}=\text{curl}(\bm{W})\text{ for some }\bm{W}\in \mathcal{V}(\bar{M})\text{ and }\text{curl}(\bm{\Gamma})=0 \}, \\
	\label{ME5}
	\mathcal{H}_N(\bar{M}):=\{\bm{\Gamma}\in \mathcal{V}(\bar{M})|\text{curl}(\bm{\Gamma})=0\text{, }\text{div}(\bm{\Gamma})=0\text{, }\bm{\Gamma}\parallel\partial\bar{M}\}.
\end{gather}
\begin{lem}
	\label{S2L2}
	Let $(\mathcal{R},g)$ be our ambient space and suppose that for every $\bar{M}\in \text{Sub}_c(\mathcal{R})\setminus\{\mathcal{R}\}$ we have $\mathcal{H}_N(\bar{M})\subset \mathcal{H}_{\text{co}}(\bar{M})$. Then $(\mathcal{R},g)$ has the div-free extension property.
\end{lem}
In Euclidean $3$-space we have $\mathcal{H}_N(\bar{M})\subset \mathcal{H}_{\text{co}}(\bar{M})$ for every compact, smoothly bounded domain $\bar{M}$, see for instance \cite[Chapter III A]{CDG01}. The following provides a broader class of such manifolds
\begin{lem}
	\label{S2C3}
	Let $\left(\mathcal{R},g\right)$ be our ambient space. If there exists a sequence $(B_n)_{n\in \mathbb{N}}\subset \mathcal{R}$ such that
	\begin{enumerate}
		\item $B_n\in \text{Sub}_c(\mathcal{R})$ for all $n\in \mathbb{N}$,
		\item $B_n\subseteq B_{n+1}$ for all $n\in \mathbb{N}$,
		\item $H^2_{dR}(B_n)=\{0\}$ for all $n\in \mathbb{N}$,
		\item $\mathcal{R}=\bigcup_{n\in \mathbb{N}}\text{int}(B_n)$,
	\end{enumerate}
	then $\mathcal{H}_N(\bar{M})\subset \mathcal{H}_{\text{co}}(\bar{M})$ for all $\bar{M}\in \text{Sub}_c(\mathcal{R})\setminus \{\mathcal{R}\}$. In particular, the ambient spaces $\mathbb{R}^3$, $\mathcal{H}^3$, $S^3$ as well as the half sphere $S^3_+$ have the div-free extension property.
\end{lem}
Here we denote by $H^2_{dR}$ the $2$-nd de Rham cohomology group. Let us also point out that the conditions in \cref{S2C3} are all of topological nature so that the divergence-free extension property in these cases holds independent of the chosen metric. We first prove \cref{S2C3}.
\newline
\newline
\underline{Proof of \cref{S2C3}:}
Fix any $\bar{M}\in \text{Sub}_c(\mathcal{R})\setminus\{\mathcal{R}\}$. By assumption on the sequence $(B_n)_n$, their interiors form an open cover of $\bar{M}$ so that by compactness of $\bar{M}$ and the inclusion property we find some $m\in \mathbb{N}$ such that $\bar{M}\subseteq B_m$. Now for any $\bm{\Gamma}\in \mathcal{H}_N(\bar{M})$ we may set it to zero outside of $\bar{M}$ to obtain a vector field $\widetilde{\bm{\Gamma}}\in L^2(B_m)$. One easily confirms that this vector field is $L^2$-orthogonal on $B_m$ to the vector space $\{\text{grad}(f)|f\in H^1_0(B_m)\}$. It then follows from the Hodge-Friedrichs decomposition \cite[Theorem 2.4.2 \& Theorem 2.4.8]{S95} that $\widetilde{\bm{\Gamma}}=\text{curl}(\bm{A})+\bm{H}$ for some $\bm{A}\in H^1(B_m)$ and $\bm{H}\in \mathcal{H}_D(B_m)$, see (\ref{S2E2}). However, according to \cite[Theorem 2.6.1 \& Corollary 2.6.2]{S95} we have $\mathcal{H}_D(B_m)\cong H^2_{dR}(B_m)=\{0\}$, the latter by assumption. We conclude $\bm{H}=0$ and so upon restricting $\bm{A}$ to $\bar{M}$ we conclude $\bm{\Gamma}\in \mathcal{H}_{\text{co}}(\bar{M})$ as desired. The last statement follows from \cref{S2L2}.
$\square$
\newline
\newline
For a detailed discussion of the div-free extension property in the Euclidean setting we refer to \cite{KMPT00}. Here we give briefly an alternative, but related, proof which applies to all abstract manifolds as described in \cref{S2L2}. First let us introduce some more notation before stating a simple lemma, needed for the proof of \cref{S2L2}
\begin{align}
	\label{S2E2}
	\mathcal{H}_D(\bar{M}):=\{\bm{\Gamma}\in \mathcal{V}(\bar{M})|\text{div}(\bm{\Gamma})=0\text{, }\text{curl}(\bm{\Gamma})=0\text{ and }\bm{\Gamma}\perp \partial\bar{M} \}, \\
	\mathcal{H}_{\text{ex}}(\bar{M}):=\{\text{grad}(f)| f\in C^{\infty}(\bar{M})\text{ and }\Delta_g f=0 \},
\end{align}
where $\Delta_g$ denotes the Laplace-Beltrami operator, i.e. $\Delta_g=-\text{div}\circ\text{grad}$ acting on functions.
\begin{lem}
	\label{S2L3}
	Let $(\bar{M},g)$ be a compact, connected, oriented, smooth Riemannian $3$-manifold with non-empty boundary and let $N:=\#\partial\bar{M}$. Then
	\[
	\dim\left(\mathcal{H}_D(\bar{M})\cap \mathcal{H}_{\text{ex}}(\bar{M}) \right)=N-1.
	\]
\end{lem}
\underline{Proof of \cref{S2L3}:} For given $1\leq i\leq N$ denote by $B_i$ the distinct boundary components of $\partial\bar{M}$. Then one can always solve the following boundary value problems, \cite[Theorem 3.4.6]{S95}, for $i=1,\dots,N-1$
\[
\Delta_g f_i=0 \text{ on }\bar{M}\text{ and }(f_i)|_{B_j}=\delta_{ij}\text{ for }1\leq j\leq N.
\]
One can then easily check that $\text{grad}(f_i)$, $i=1,\dots,N-1$ are all linearly independent and in fact form a basis of the space of interest. $\square$
\newline
\newline
Now we prove \cref{S2L2}
\newline
\newline
\underline{Proof of \cref{S2L2}:} We may consider the flowout of $\partial\bar{M}$ along the outward pointing unit normal, see also \cite[Theorem 9.24]{L12}, and obtain for small enough times some open $\bar{M}\subset U\subseteq \mathcal{R}$, such that $\overline{U}\in \text{Sub}_c(\mathcal{R})$ (here $\overline{U}$ denotes the topological closure of $U$ in $\mathcal{R}$) and such that $\overline{U}$ deformation retracts onto $\bar{M}$. Since the flowout provides for every fixed time a diffeomorphism between $\partial\bar{M}$ and its images, we see that in particular $\overline{U}$ and $\bar{M}$ have the same number of connected components. Let as before $N:=\#\partial\bar{M}$, then we can fix a basis $\bm{\Gamma}_i\in \mathcal{H}_D(\overline{U})\cap \mathcal{H}_{\text{ex}}(\overline{U})$, $i=1,\dots, N-1$ and denote by $\gamma_i$ the associated $1$-forms via the metric $g$. Further, let $\iota$ denote the inclusion map between the boundary of a manifold with boundary and the manifold itself and by $\iota^\#$ its pullback. Then we may consider the $2$-forms $\iota^{\#}\star\gamma_i\in \Omega^2(\partial\overline{U})$, where $\Omega^k$ denotes the space of smooth $k$-forms on a given manifold and $\star$ denotes the Hodge star operator. We now first observe that if for some $\bm{\Gamma}\in \mathcal{H}_D(\overline{U})$, the associated $2$-form $\iota^{\#}\star \gamma\in \Omega^{2}(\partial\overline{U})$ integrates to zero on every connected component of $\partial\overline{U}$, then it is exact on the boundary, i.e. there exists some $\alpha\in \Omega^{1}(\partial\overline{U})$ with $\iota^{\#}\star \gamma=d\alpha$, \cite[Chapter 8 Theorem 9]{Sp99}. It then follows from the Hodge decomposition theorem \cite[Corollary 3.5.2]{S95} and our assumption $\mathcal{H}_N(\overline{U})\subset \mathcal{H}_{\text{co}}(\overline{U})$ that $\mathcal{H}_D(\overline{U})\subset \mathcal{H}_{\text{ex}}(\overline{U})$, which in turn implies that the associated $2$-form $\star \gamma \in \Omega^{2}(\overline{U})$ is coexact, i.e. there is some $\beta\in \Omega^3(\overline{U})$ with $\star\gamma=\delta \beta$. We then conclude from the integration by parts formula and by extending the form $\alpha\in \Omega^{1}(\partial\overline{U})$ to some smooth form $\tilde{\alpha}\in \Omega^{1}(\overline{U})$ with $\iota^{\#}\tilde{\alpha}=\alpha$, that $\bm{\Gamma}$ must be already the zero vector field, whenever its associated $2$-form integrates to zero on each boundary component of $\overline{U}$.
\newline
Coming back to our linearly independent $2$-forms $\iota^{\#}\star\gamma_i$, $i=1,\dots,N-1$, the preceding considerations imply that the matrix $C\in \mathbb{R}^{(N-1)\times (N-1)}$, whose entries are given by $C_{ji}:=\int_{\partial U_j}\iota^{\#}\star\gamma_i$ ($\partial U_j$ denoting the $j$-th boundary component of $\overline{U}$), is invertible. Further, since the flowout gives rise to a smoothly varying family of diffeomorphisms, connecting $\partial\bar{M}$ and $\partial\overline{U}$, we can conclude that $\int_{\partial\bar{M}_j}\iota^{\#}\star\gamma_i=\int_{\partial\overline{U}_j}\iota^{\#}\star\gamma_i=C_{ji}$, $1\leq j,i\leq N-1$, where $\partial\bar{M}_j$ denote the corresponding boundary components of $\bar{M}$ diffeomorphic to $\partial\overline{U}_j$ by means of the flowout and $\iota$ in the first integral denotes the inclusion $\iota:\partial\bar{M}\rightarrow \overline{U}$.
\newline
Now given any divergence-free $\bm{Y}\in \mathcal{V}(\bar{M})$, we define $a_j:=\int_{\partial\bar{M}_j}\iota^{\#}\omega^2_{\bm{Y}}$, where $\omega^2_{\bm{Y}}\in \Omega^2(\bar{M})$ denotes the corresponding $2$-form and we set $\vec{a}:=(a_1,\dots,a_{N-1})$. By invertibility of the matrix $C$ we may find some vector $\vec{\lambda}\in \mathbb{R}^{N-1}$ with $C\cdot\vec{\lambda}=\vec{a}$, whose entries are denoted by $\lambda_i$. Define the vector field
\[
\bm{X}:=\bm{Y}-\sum_{i=1}^{N-1}\lambda_i\bm{\Gamma}_i.
\]
By construction we see that the associated $2$-form $\iota^{\#}\omega^2_{\bm{X}}$ integrates to zero on every boundary component $\partial\bar{M}_j$ for $1\leq j\leq N-1$. Since $\bm{X}$ is divergence-free it follows from Stokes-theorem that the corresponding $2$-form must also integrate to zero on the remaining boundary component. It then follows from the Hodge-decomposition theorem \cite[Corollary 2.4.9]{S95} that there exists some $\bm{\Gamma}\in \mathcal{H}_D(\bar{M})$ with associated $2$-form $\star\gamma$, such that $\omega^2_{\bm{X}}=d\alpha+\star\gamma$ for some suitable $\alpha\in \Omega^1(\bar{M})$. Note that every exact $2$-form integrates to zero on every boundary component, so that we conclude that $\star \gamma$ must integrate to zero on every boundary component, since $\omega^2_{\bm{X}}$ does. An identical argument as in the case of $\overline{U}$ implies that $\star\gamma=0$ and consequently $\omega^2_{\bm{X}}=d\alpha$. Extending $\alpha$ to a smooth form on $U$, which is always possible by means of Seeley's extension theorem \cite{See64}, we see that $\bm{X}$ admits a divergence-free extension to $U$ and therefore, since the $\bm{\Gamma}_i$ are defined on $U$, so does $\bm{Y}$. $\square$
\newline
\newline
We now come to the notion of concircular vector fields. The following definition is standard
\begin{defn}[Concircular vector fields]
	\label{S2D5}
	Let $(\mathcal{R},g)$ be an oriented, smooth Riemannian $3$-manifold with (possibly empty) boundary. We say that a vector field $\bm{Z}\in \mathcal{V}(\mathcal{R})$ is \textit{concircular} if there exists a smooth function $f\in C^{\infty}(\mathcal{R})$, the so called potential function, such that for all vector fields $\bm{X}\in \mathcal{V}(\mathcal{R})$ we have
	\[
	\nabla_{\bm{X}}\bm{Z}=f\bm{X},
	\]
	where $\nabla$ denotes the Levi-Civita connection.
\end{defn}
We will give some examples after proving the main result of this section. First we will need the following result, which follows from the arguments of the proofs of \cite[Theorem B]{CDGT002} and \cite[Theorem 2.2.3 and Proposition 2.2.4]{G20Diss} in combination with the Hodge-decomposition theorem \cite[Corollary 2.4.9]{S95}
\begin{lem}
	\label{S2L6}
	Let $(\mathcal{R},g)$ be our ambient space and assume that it has the div-free extension property. If $\bar{M}\in \text{Sub}^V_c(\mathcal{R})$ for some $0<V<\text{vol}(\mathcal{R})$ is optimal, then for every $\bm{X}\in \mathcal{V}_n(\bar{M})$ with $\text{curl}(\bm{X})=\lambda_+(\bar{M})\bm{X}$ we have
	\[
	|\bm{X}|_g=c\text{ on }\partial\bar{M}\text{ for some constant }c\geq 0.
	\]
\end{lem}
Note that the constant $c$ has the same value on every connected component of the boundary. The first main observation of this section is the following, which establishes a connection between the existence of concircular vector fields and the value $c$ on the boundary.
\begin{lem}
	\label{S2L7}
	Let $(\mathcal{R},g)$ be our ambient space. Assume that $(\mathcal{R},g)$ has the div-free extension property and that there exists some concircular vector field $\bm{Z}\in \mathcal{V}(\mathcal{R})$ with potential function $f$. If $\bar{M}\in \text{Sub}^V_c(\mathcal{R})$ is an optimal domain for some $0<V<\text{vol}(\mathcal{R})$ and $\bm{X}\in \mathcal{V}_n(\bar{M})$ satisfies $\text{curl}(\bm{X})=\lambda_+(\bar{M})\bm{X}$, then the constant $c^2=g|_{\partial\bar{M}}(\bm{X},\bm{X})$ satisfies
	\[
	3c^2\int_{\bar{M}}f\omega_g=\int_{\bar{M}}f g(\bm{X},\bm{X})\omega_g,
	\]
	where $\omega_g$ denotes the volume form induced by the metric $g$.
\end{lem}
\underline{Proof of \cref{S2L7}:} We start by computing
\[
2\langle \bm{X},f\bm{X}\rangle_{L^2(\bar{M})}=2\langle \bm{X},\nabla_{\bm{X}}\bm{Z}\rangle_{L^2(\bar{M})}=-2\langle \nabla_{\bm{X}}\bm{X},\bm{Z}\rangle_{L^2(\bar{M})},
\]
where we used an integration by parts formula and the fact that $\bm{X}$ is divergence-free and tangent to the boundary. We can now use the following calculus formula
\[
2\nabla_{\bm{X}}\bm{X}=\text{grad}(g(\bm{X},\bm{X}))-2\bm{X}\times \text{curl}(\bm{X})=\text{grad}(g(\bm{X},\bm{X})),
\]
where the last identity follows because $\bm{X}$ is a strong Beltrami field. We conclude by means of integration by parts and the Stokes' theorem
\[
2\langle \bm{X},f\bm{X}\rangle_{L^2(\bar{M})}=-\langle \text{grad}(g(\bm{X},\bm{X})),\bm{Z}\rangle_{L^2}=\langle g(\bm{X},\bm{X}),\text{div}(\bm{Z})\rangle_{L^2}-\int_{\partial\bar{M}}g(\bm{X},\bm{X})\iota^{\#}\omega^2_{\bm{Z}}
\]
\[
=\langle g(\bm{X},\bm{X}),\text{div}(\bm{Z})\rangle_{L^2}-c^2\int_{\bar{M}}\text{div}(\bm{Z})\omega_g=3\int_{\bar{M}}fg(\bm{X},\bm{X})\omega_g-3c^2\int_{\bar{M}}f\omega_g,
\]
where we used that $g(\bm{X},\bm{X})=c^2$ on the boundary and that $\text{div}(\bm{Z})=3f$ (since we are in the $3$-dimensional setting and $\bm{Z}$ is concircular). $\square$
\subsection{Preparing for the proofs}
In the upcoming proposition we denote for given $\bar{M}$ by $|\partial\bar{M}_j|$ the surface area of the respective boundary components and by a signed sum of the boundary components we mean any sum of the form $\sum_{j=1}^{\#\partial\bar{M}}\pm|\partial\bar{M}_j|$, where the sign might differ from summand to summand. Also recall that $E_+(\bar{M})$ denotes the eigenspace corresponding to the eigenvalue $\lambda_+(\bar{M})$. With this at hand we can prove the following.
\begin{prop}
	\label{S2P8}
	Let $(\mathcal{R},g)$ be our ambient space. Assume that $(\mathcal{R},g)$ has the div-free extension property and that it admits a concircular vector field $\bm{Z}\in \mathcal{V}(\mathcal{R})$ with potential function $f$. Suppose further that $\bar{M}\in \text{Sub}^V_c(\mathcal{R})$ is an optimal domain for some $0<V<\text{vol}(\mathcal{R})$ and that $\int_{\bar{M}}f\omega_g\neq 0$. Then the following holds
	\begin{enumerate}
	\item $1\leq \dim(E_+(\bar{M}))\leq 2$.
	\item If no signed sum of the boundary components of $\bar{M}$ is zero, then $\lambda_+(\bar{M})$ is simple.
	\item If $\dim(E_+(\bar{M}))= 2$, then $\left(\partial\bar{M},\iota^\#g\right)$ is flat, where $\iota:\partial\bar{M}\rightarrow\mathcal{R}$ as usual denotes the pullback metric.
\end{enumerate}
\end{prop}
\underline{Proof of \cref{S2P8}:} We assume that $\int_{\bar{M}}f\omega_g\neq 0$ and show that in this case the eigenvalue has at most multiplicity $2$, respectively $1$ under the additional no-zero-signed sum assumption. To simplify notation we define $I_f(\bm{V},\bm{W}):=\int_{\bar{M}}fg(\bm{V},\bm{W})\omega_g$ for any two given smooth vector fields $\bm{V}$ and $\bm{W}$ on $\bar{M}$. Now suppose for the moment that there exist two linearly independent curl eigenfields $\bm{X}_1,\bm{X}_2\in \mathcal{V}_n(\bar{M})$ corresponding to the eigenvalue $\lambda_+(\bar{M})$. We define $c^2_i:=g|_{\partial\bar{M}}(\bm{X}_i,\bm{X}_i)$, $i=1,2$ and observe that these constants must be non-zero, since the $\bm{X}_i$ are linearly independent, see \cite[Vainshtein's lemma]{CDGT002} and \cite[Lemma 2.1]{G21a}. By \cref{S2L7} we therefore have $I_f(\bm{X}_i,\bm{X}_i)\neq 0$ for $i=1,2$ and after an appropriate rescaling we can always achieve $I_f(\bm{X}_i,\bm{X}_i)=3\int_{\bar{M}}f\omega_g$, so that in turn by \cref{S2L7} we have $c^2_i=1$ for $i=1,2$. Now a Gram-Schmidt argument implies that we can always force $I_f(\bm{X}_1,\bm{X}_2)=0$. We then define the following vector field $\bm{X}:=\frac{1}{\sqrt{2}}\left(\bm{X}_1+\bm{X}_2\right)$ and observe that it is yet again an eigenfield of curl corresponding to the eigenvalue $\lambda_+(\bar{M})$. Further we have
\[
I_f(\bm{X},\bm{X})=\frac{I_f(\bm{X}_1,\bm{X}_1)+I_f(\bm{X}_2,\bm{X}_2)}{2}=3\int_{\bar{M}}f\omega_g,
\]
so that by means of \cref{S2L7} we accordingly have
\[
1=c^2\equiv g|_{\partial\bar{M}}(\bm{X},\bm{X})=\frac{g|_{\partial\bar{M}}(\bm{X}_1,\bm{X}_1)+g|_{\partial\bar{M}}(\bm{X}_2,\bm{X}_2)}{2}+g|_{\partial\bar{M}}(\bm{X}_1,\bm{X}_2)
\]
\[
=\frac{c^2_1+c^2_2}{2}+g|_{\partial\bar{M}}(\bm{X}_1,\bm{X}_2)=1+g|_{\partial\bar{M}}(\bm{X}_1,\bm{X}_2)
\]
and consequently $g|_{\partial\bar{M}}(\bm{X}_1,\bm{X}_2)=0$ on $\partial\bar{M}$. Since $g|_{\partial\bar{M}}(\bm{X}_i,\bm{X}_i)=c^2_i=1$ on the boundary, we conclude that the vector fields $\bm{X}_i$, $i=1,2$ are everywhere linearly independent on $\partial\bar{M}$. Now if there exist three linearly independent eigenfields $\bm{X}_i$, $i=1,2,3$, then just like before we can apply the Gram-Schmidt procedure to ensure that $I_f(\bm{X}_i,\bm{X}_j)=0$ for $i\neq j$ and after scaling $c^2_i=1$ for $i=1,2,3$. Then an identical argument as above implies $g|_{\partial\bar{M}}(\bm{X}_i,\bm{X}_j)=0$ for $i\neq j$ and consequently, since all the $\bm{X}_i$ are tangent to the boundary, we found $3$-linearly independent vector fields tangent to the $2$-dimensional boundary, which is a contradiction. Hence the eigenspace is at most $2$-dimensional.
\newline
\newline
From now on we additionally assume that the no-zero-signed sum condition is imposed and that our eigenspace is $2$-dimensional. As already shown we can find eigenfields $\bm{X}_i$, $i=1,2$, which are $g$-orthogonal and of unit speed on the boundary. We observe that by means of standard calculus formulae we have
\[
\text{div}(\bm{X}_1\times \bm{X}_2)=g(\text{curl}(\bm{X}_1),\bm{X}_2)-g(\bm{X}_1,\text{curl}(\bm{X}_2))=\lambda_+(\bar{M})\left(g(\bm{X}_1,\bm{X}_2)-g(\bm{X}_1,\bm{X}_2) \right)=0.
\]
Further, since $\bm{X}_1$ and $\bm{X}_2$ are everywhere $g$-orthogonal and tangent to the boundary, we must have
\[
\bm{X}_1\times \bm{X}_2=\pm |\bm{X}_1\times \bm{X}_2|_g\mathcal{N}=\pm \mathcal{N},
\]
on every boundary component, where $\mathcal{N}$ denotes the outward pointing unit normal and where we used that $|\bm{X}_i|_g=1$ for $i=1,2$. So if we let $N:=\#\partial\bar{M}$, we find
\[
0=\int_{\bar{M}}\text{div}(\bm{X}_1\times \bm{X}_2)\omega_g=\int_{\partial\bar{M}}g(\bm{X}_1\times \bm{X}_2,\mathcal{N})\omega_{g_{\partial\bar{M}}}=\sum_{j=1}^N\int_{\partial\bar{M}_j}g(\bm{X}_1\times \bm{X}_2,\mathcal{N})\omega_{g_{\partial\bar{M}}}
\]
\[
=\sum_{j=1}^N\pm \int_{\partial\bar{M}_j}g(\mathcal{N},\mathcal{N})\omega_{g_{\partial\bar{M}}}=\sum_{j=1}^N\pm |\partial\bar{M}_j|,
\]
where $\omega_{g_{\partial\bar{M}}}$ denotes the volume form on $\partial\bar{M}$ induced by $g$ and where $|\partial\bar{M}_j|$ denotes the surface area of the corresponding boundary component. Here the sign $\pm$ may differ from component to component. This contradicts our assumption and so we conclude that the eigenvalue $\lambda_+(\bar{M})$ must be simple.
\newline
\newline
As for the last claim, let us suppose that $\dim(E_+(\bar{M}))=2$. We can as before find two eigenfields $\bm{X}_1,\bm{X}_2\in E_+(\bar{M})$ which are of unit speed and $g$-orthogonal at the boundary. We can then use the following vector calculus identity
\[
\text{grad}(g(\bm{X},\bm{Y}))=\nabla_{\bm{X}}\bm{Y}+\nabla_{\bm{Y}}\bm{X}+\bm{X}\times \text{curl}(\bm{Y})+\bm{Y}\times \text{curl}(\bm{X})\text{ for all }\bm{X},\bm{Y}\in \mathcal{V}(\bar{M}).
\]
Setting $\bm{X}=\bm{X}_1$, $\bm{Y}=\bm{X}_2$ and using the eigenfield property we find $\text{grad}(g(\bm{X}_1,\bm{X}_2))=\nabla_{\bm{X}_1}\bm{X}_2+\nabla_{\bm{X}_2}{\bm{X}_1}$. Using that $g(\bm{X}_1,\bm{X}_2)=0$ on $\partial\bar{M}$ we obtain
\[
\left(\nabla_{\bm{X}_1}\bm{X}_2\right)^\parallel=-\left(\nabla_{\bm{X}_2}{\bm{X}_1}\right)^\parallel\text{ on }\partial\bar{M},
\]
where $\bm{X}^\parallel$ denotes the tangential part of a given vector field $\bm{X}\in \mathcal{V}(\bar{M})$ on $\partial\bar{M}$. Now by means of the defining properties of the Levi-Civita connection we have the identity, recall that $\bm{X}_1$ is tangent to the boundary,
\[
0=g\left(\bm{X}_1,\text{grad}(g(\bm{X}_1,\bm{X}_2))\right)=g\left(\nabla_{\bm{X}_1}\bm{X}_1,\bm{X}_2\right)+g\left(\bm{X}_1,\nabla_{\bm{X}_1}\bm{X}_2\right)=g\left(\bm{X}_1,\nabla_{\bm{X}_1}\bm{X}_2\right)\text{ on }\partial\bar{M},
\]
where we used that $2\nabla_{\bm{X}_1}\bm{X}_1=\text{grad}(|\bm{X}_1|^2)$ by setting $\bm{X}=\bm{X}_1=\bm{Y}$ in the above identity, that $|\bm{X}_1|=1$ on the boundary and that $\bm{X}_2$ is tangent to the boundary. Reversing the roles of $\bm{X}_1$ and $\bm{X}_2$, using the tangency and the derived relation between $\nabla_{\bm{X}_1}\bm{X}_2$ and $\nabla_{\bm{X}_2}\bm{X}_1$ we obtain
\[
0=g(\bm{X}_2,\nabla_{\bm{X}_2}\bm{X}_1)=g\left(\bm{X}_2,\left(\nabla_{\bm{X}_2}\bm{X}_1\right)^\parallel\right)=-g\left(\bm{X}_2,\left(\nabla_{\bm{X}_1}\bm{X}_2\right)^\parallel\right)=-g\left(\bm{X}_2,\nabla_{\bm{X}_1}\bm{X}_2\right).
\]
Since $\bm{X}_1,\bm{X}_2$ span $T_p\partial\bar{M}$ at each $p\in \partial\bar{M}$ we conclude that $\nabla_{\bm{X}_1}\bm{X}_2=f_1\mathcal{N}$ on $\partial\bar{M}$ for a suitable smooth function $f_1\in C^\infty(\partial\bar{M})$ and $\mathcal{N}$ is the outward pointing unit normal. In the same spirit we find a function $f_2\in C^\infty(\partial\bar{M})$ with $\nabla_{\bm{X}_2}\bm{X}_1=f_2\mathcal{N}$ and so letting $f:=f_1-f_2$ we find
\[
[\bm{X}_1,\bm{X}_2]=f\mathcal{N}\perp\partial\bar{M},
\]
where $[\cdot,\cdot]$ denotes the Lie-bracket of vector fields. But it is a standard fact that the Lie-bracket of two vector fields which are tangent to the boundary is itself tangent to the boundary. Hence we conclude $[\bm{X}_1,\bm{X}_2]=0$ on $\partial\bar{M}$. Thus, upon viewing the vector fields $\bm{X}_i$ as vector fields on $\partial\bar{M}$ for $i=1,2$, we found (with a slight abuse of notation) two vector fields $\bm{X}_i$, $i=1,2$, on $\partial\bar{M}$ which are everywhere linearly independent and which commute with each other. Hence we can find around any given $p\in \partial \bar{M}$ a coordinate chart $\mu_p$ of $\partial\bar{M}$ around $p$ such that the basis vectors satisfy $\partial_i=\bm{X}_i$ for $i=1,2$ on the domain of the chart. Consequently, as $g(\bm{X}_i,\bm{X}_j)=\delta_{ij}$ on $\partial\bar{M}$, where $\delta_{ij}$ is the Kronecker delta, we see that the metric tensor $\iota^\#g$ is locally given by $g_{ij}=\delta_{ij}$ in these coordinate charts. We conclude that $\left(\partial\bar{M},\iota^\#g\right)$ is flat as claimed. $\square$
\newline
\newline
Let us shortly remark that given any compact domain $\bar{M}$ with boundary components $\partial\bar{M}_j$, $j=1,\dots,N$ and if there exists some $1\leq i \leq N$ with $|\partial\bar{M}_i|>\sum_{j\neq i}|\partial\bar{M}_j|$, then necessarily no signed sum of the boundary components can be zero.
\subsection{The proofs}
\underline{Proof of \cref{MT2}:} By \cref{S2C3} and \cref{S2P8} it is enough to show the existence of concircular vector fields with strictly positive potential functions and that no closed (orientable), flat surface embeds isometrically into the corresponding ambient spaces. The fact that the Gaussian curvature of each closed, embedded surface in $\mathbb{R}^3$, $\mathcal{H}^3$ or $S^3_+$ must be somewhere non-zero is standard and for the convenience of the reader we provide an argument in the appendix, see \cref{AL1}. Hence it is enough to provide concircular vector fields with positive potential functions.
\newline
\newline
(i) The case $\mathbb{R}^3$: The position vector field $\vec{x}$ is a concircular vector field with potential function $f\equiv 1$, so that the result follows immediately from \cref{S2P8} and \cref{AL1}.
\newline
\newline
(ii) The case $\mathcal{H}^3$: We consider the Poincar\'{e} ball model and once again the position vector field. One readily checks that it is a concircular vector field with potential function $f(\vec{x})=\frac{1+|\vec{x}|^2_{2}}{1-|\vec{x}|^2_{2}}$, where $|\cdot|_{2}$ denotes the Euclidean distance and we recall that we identify the hyperbolic space with the open unit ball. Hence the potential function is strictly positive and \cref{S2P8} in combination with \cref{AL1} apply.
\newline
\newline
(iii) The case $S^3_+$: By \cref{MR6} we know that $S^3_+$ is isometric to the open unit ball with metric $g(\vec{x})=\frac{4}{\left(1+|\vec{x}|^2_{2}\right)^2}g_E(\vec{x})$. Once again we check that the position vector field is concircular with potential function $f(\vec{x})=\frac{1-|\vec{x}|^2_{2}}{1+|\vec{x}|^2_{2}}$, which is strictly positive on the open unit ball, so that once more \cref{S2P8} and \cref{AL1} apply. This concludes the proof of \cref{MT2}. $\square$
\newline
\newline
Let us point out that Clifford tori provide examples of flat surfaces embedded in $S^3$, so that we have to make the additional assumption of a connected boundary in the statement of \cref{MT3}.
\newline
\newline
\underline{Proof of \cref{MT3}:} Again \cref{S2C3} implies that $S^3$ has the div-free extension property. Since we assume $0<V<\text{vol}(S^3)$, after applying an isometry, we can always achieve that the north pole is not contained in $\bar{M}$. Then we may identify $\bar{M}$ with a compact domain in $\mathbb{R}^3$ by means of the stereographic projection, where $\mathbb{R}^3$ is equipped with the metric $g(\vec{x})=\frac{4}{\left(1+|\vec{x}|^2_{2}\right)^2}g_E(\vec{x})$. We note that the optimality of $\bar{M}$ on $S^3$ in particular implies the optimality of $\bar{M}$ in $\mathbb{R}^3$ (with the mentioned metric). We already know from the proof of \cref{MT2} that the position vector field is a concircular vector field with potential function $f(\vec{x})=\frac{1-|\vec{x}|^2_{2}}{1+|\vec{x}|^2_{2}}$. Now if $\lambda_+(\bar{M})$ is not simple it follows from \cref{S2P8} that the function $f$ must integrate to zero over the optimal domain. Pulling back the function $f$ by means of the stereographic projection yields the result. $\square$
\section{Proof of \cref{MT5}}
\Cref{MT5} will follow from the following more general result, which also includes the Euclidean case as a corollary. Here we denote by $H^k_{dR}(\bar{M})$ the $k$-th de Rham cohomology group of the manifold $\bar{M}$.
\begin{prop}
	\label{S5P1}
	Let $(\mathcal{R},g)$ be our ambient space and assume that $(\mathcal{R},g)$ has the div-free extension property and admits a concircular vector field with strictly positive potential function. Suppose further that $(\mathcal{R},g)$ admits a Killing field $\bm{Y}\in \mathcal{V}(\mathcal{R})$ (i.e. it satisfies the Killing equations) with the following three properties
	\begin{enumerate}
		\item $g(\bm{Y},\text{curl}(\bm{Y}))=0$,
		\item $\forall c^2\geq 0$ the set $\{p\in \mathcal{R}||\bm{Y}|^2_g=c^2 \}$ does not contain a subset homeomorphic to $T^2=S^1\times S^1$,
		\item $\{p\in \mathcal{R}|\bm{Y}(p)=0 \}$ does not contain a subset homeomorphic to $S^1$.
	\end{enumerate}
Given some $0<V<\text{vol}(\mathcal{R})$, if there exists an optimal domain $\bar{M}\in \text{Sub}^V_c(\mathcal{R})$ with $H^2_{dR}(\bar{M})=\{0\}$ and such that $\bm{Y}\parallel \partial\bar{M}$, then the set
\begin{align}
	\label{S5E1}
\partial\bar{M}\setminus \overbrace{\{p\in \partial\bar{M}| |\bm{Y}(p)|_g=\min_{q\in \partial\bar{M}}|\bm{Y}(q)|_g \}}^{=:N_0}
\end{align}
is disconnected and $\bm{Y}(p)\neq 0$ for every $p\in \bar{M}$. In addition, if $\bm{X}\in \mathcal{V}_n(\bar{M})$ is any fixed eigenfield of $\lambda_+(\bar{M})$, then $\bm{X}$ and $\bm{Y}$ are linearly dependent on $\partial\bar{M}$ precisely on the set $N_0$.
\end{prop}
Note that the proof of \cref{MT2} showed that the spaces $\mathbb{R}^3$ and $\mathcal{H}^3$ admit concircular vector fields with strictly positive potential functions and that \cref{S2C3} shows that these spaces have the div-free extension property. It is then not hard to verify that the vector field $\bm{R}$ defined in (\ref{ME4}) satisfies the requirements (i)-(iii) of \cref{S5P1} for both of these spaces and their respective metrics. Further, it follows from \cite{CDG02} that the connectedness of the boundary implies $H^2_{dR}(\bar{M})=\{0\}$. Lastly, the proof of \cref{S5P1} shows that $\bm{Y}(p)\neq 0$ on $\bar{M}$ for any not-identically zero Killing field $\bm{Y}$ with properties (i)-(iii) and any optimal domain (with not necessarily connected boundary), as long as $\bm{Y}\parallel \partial\bar{M}$. Therefore \cref{MT5} as well as the corresponding Euclidean version follow immediately from \cref{S5P1}.
\begin{rem}
	\label{S5R2}
\begin{enumerate}
\item It follows from assumption (ii) of \cref{S5P1} that the Killing field $\bm{Y}$ is not the zero vector field. Further, it follows from \cite[Theorem 8.1.5]{Pe16} that the zero set is a $1$-dimensional manifold. In this sense condition (iii) is not as restrictive as it appears on first glance.
\item For every connected, compact, oriented $3$-manifold with non-empty boundary, we have the implication: If $H^2_{dR}(\bar{M})=0$, then $\partial\bar{M}$ is connected, which follows for example from \cite[Lemma 2.2]{G21a}. Further, if we do not demand that $(\mathcal{R},g)$ has the div-free extension property, then \cref{S2L6} still holds, with the only difference that the constant $c$ may differ from boundary component to boundary component. However, since our assumptions imply that the boundary is connected, we may in fact drop the div-free extension property assumption in \cref{S5P1} and obtain the exact same result. The first part of the proof, during which we do not use the assumption $H^2_{dR}(\bar{M})=\{0\}$, is still valid, even if the constant has possibly distinct values on different boundary components.
	\end{enumerate}
\end{rem}
\underline{Proof of \cref{S5P1}:} First we note that $\partial\bar{M}\neq \emptyset$. Otherwise $\text{int}(\bar{M})=\bar{M}$ implying that $\bar{M}$ is an open subset of $\mathcal{R}$ (recall that topological and manifold interiors\slash boundaries coincide for all $\bar{M}\in \text{Sub}_c(\mathcal{R})$).Then the compactness implies that $\bar{M}$ is an open and closed subset of $\mathcal{R}$ and the connectedness of $\mathcal{R}$ yields $\bar{M}=\mathcal{R}$, which contradicts the fact that $\text{vol}(\bar{M})=V<\text{vol}(\mathcal{R})$. Thus, indeed $\partial\bar{M}\neq \emptyset$. We will now argue that $\bm{Y}(p)\neq 0$ on $\partial\bar{M}$. To this end we note that, since $\bm{Y}\parallel \partial\bar{M}$, we may view it as a vector field on $\partial\bar{M}$ which again is Killing with respect to the induced metric. It then follows that either $\bm{Y}$ has only isolated zeros on $\partial\bar{M}$ or vanishes identically on at least one boundary component of $\partial\bar{M}$ \cite[Proposition 8.1.4 \& Theorem 8.1.5]{Pe16}. Now if $\bm{Y}$ vanishes identically on a connected component of $\partial\bar{M}$, then once more \cite[Proposition 8.1.4 \& Theorem 8.1.5]{Pe16} imply, since $\bm{Y}$ is Killing on the ambient space, that $\bm{Y}$ vanishes on all of $\mathcal{R}$ and consequently that $\bm{Y}$ is the zero vector field, which contradicts assumption (ii) of \cref{S5P1}. Now assume that $\bm{Y}$ has a zero at some boundary component of $\partial\bar{M}$, then since all these zeros are isolated and since all isolated zeros of Killing fields have index $+1$, it follows from the Poincar\'{e}-Hopf theorem, \cite[Chapter 6]{M65}, and the classification of closed (orientable) surfaces, \cite[Chapter 9 Theorem 3.5]{H76}, that the considered component must be a sphere. Now \cref{S2L6}, the fact that any eigenfield $\bm{X}\in \mathcal{V}_n(\bar{M})\setminus\{0\}$ of $\lambda_+(\bar{M})$ is tangent to the boundary and the fact that a boundary component is a sphere, implies that we must have $\bm{X}\equiv 0$ on $\partial\bar{M}$. Now a unique continuation result for curl eigenfields, \cite[Lemma 2.1]{G21a} implies that $\bm{X}\equiv 0$ on $\bar{M}$, which is a contradiction. Hence $\bm{Y}(p)\neq 0$ on $\partial\bar{M}$. Now suppose that $\bm{Y}(p)=0$ for some $p\in M:=\text{int}(\bar{M})$. Then, since the restriction of $\bm{Y}$ to $M$ is again Killing, its zero set must be a $1$-dimensional submanifold (without boundary), \cite[Proposition 8.1.4 \& Theorem 8.1.5]{Pe16}. Now the set $\{p\in M|\bm{Y}(p)=0\}$ must be closed in $\mathcal{R}$ since $\partial\bar{M}$ does not contain any zeros. Since $\bar{M}$ is compact, so must be $\{p\in M|\bm{Y}(p)=0\}$ which shows that its connected components are all closed circles, which contradicts assumption (iii) of \cref{S5P1}. Overall $\bm{Y}\neq 0$ on $\bar{M}$.
\newline
\newline
With this at hand we may define the vector field $\bm{\Gamma}:=\frac{\bm{Y}}{|\bm{Y}|^2_g}$ on $\bar{M}$. It now follows from the Killing equations and the assumption $g(\bm{Y},\text{curl}(\bm{Y}))=0$ that $\bm{\Gamma}$ is divergence- and curl-free. Since $\bm{Y}\parallel \partial\bar{M}$, we find $\bm{\Gamma}\in \mathcal{H}_N(\bar{M})$, recall (\ref{ME5}). Hence, if we let $\bm{X}$ be any fixed eigenfield of $\lambda_+(\bar{M})$, we compute by means of Stokes' theorem
\begin{align}
	\label{S5E2}
	\int_{\partial\bar{M}}g(\bm{X}\times \bm{\Gamma},\mathcal{N})\omega_{g_{\partial\bar{M}}}=\int_{\bar{M}}\text{div}(\bm{X}\times \bm{\Gamma})\omega_g	=\lambda_+(\bar{M})\langle \bm{X},\bm{\Gamma}\rangle_{L^2}=0,
\end{align}
where we used that the space $\mathcal{V}_n(\bar{M})$ is $L^2$-orthogonal to $\mathcal{H}_N(\bar{M})$ and standard calculus identities.
\newline
\newline
From now on we assume that $H^2_{dR}(\bar{M})=\{0\}$. As pointed out in \cref{S5R2}[(ii)] this implies that $\partial\bar{M}$ is connected and since $\bm{X}$ as well as $\bm{\Gamma}$ are tangent to the boundary we see that $\bm{X}\times \bm{\Gamma}=\pm|\bm{X}\times \bm{\Gamma}|_g\mathcal{N}$, where the sign may differ at different boundary points. We conclude
\[
\int_{\partial\bar{M}}\pm|\bm{X}\times \bm{\Gamma}|\omega_{g_{\partial\bar{M}}}=0,
\]
from which it immediately follows that the set $L:=\{p\in \partial\bar{M}| \bm{X}(p)\times \bm{Y}(p)=0\}=\{p\in \partial\bar{M}| \bm{X}(p)\times \bm{\Gamma}(p)=0\}$ must either disconnect the boundary, i.e. $\partial\bar{M}\setminus L$ is not connected or otherwise $L=\partial\bar{M}$. We will now show that in fact $L=N_0$, see (\ref{S5E1}).
\newline
We have already argued that $\bm{Y}$ does not vanish on $\partial\bar{M}$. Thus, the classification of surfaces and the Poincar\'{e}-Hopf theorem imply that all boundary components of $\bar{M}$ must be tori. Hence, once we show that $L=N_0$, it will follow from assumption (ii) on $\bm{Y}$ that $L\neq \partial\bar{M}$ and the proof will be concluded.
\newline
To see this identity we recall that $\partial\bar{M}$ is connected, so that our assumptions in combination with \cref{S2P8} imply that the eigenvalue $\lambda_+(\bar{M})$ is simple. Now it follows easily from the fact that $\bm{Y}$ induces isometries that the vector fields $\bm{X}$ and $\bm{Y}$ commute. Then the Killing equations and the eigenfield property imply the following identities
\[
\text{grad}(g(\bm{X},\bm{Y}))=[\bm{Y},\bm{X}]+\bm{Y}\times \text{curl}(\bm{X})=\lambda_+(\bar{M})\bm{Y}\times \bm{X}.
\]
Since $\bm{Y}$ and $\bm{X}$ are tangent to the boundary, we see that $\text{grad}(g(\bm{X},\bm{Y}))\perp \partial\bar{M}$ and by connectedness of the boundary we must have $g|_{\partial\bar{M}}(\bm{Y},\bm{X})=c_0=\text{const}$. We have seen already that $\emptyset\neq L$ and hence that $\bm{X}$ and $\bm{Y}$ are linearly dependent in at least one point. Fix any $p\in L$, then we must have by \cref{S2L6} (after possibly rescaling $\bm{X}$ and replacing $\bm{X}$ by $-\bm{X}$)
\[
c_0=g(\bm{X}(p),\bm{Y}(p))=|\bm{X}(p)|_g|\bm{Y}(p)|_g=|\bm{Y}(p)|_g\text{ for all }p\in L.
\]
Now let $d_0:=\min_{q\in \partial\bar{M}}|\bm{Y}(q)|_g$ and suppose that $p\in L\setminus N_0$, then we find for any $q\in N_0\neq \emptyset$
\[
d_0<|\bm{Y}(p)|_g=|c_0|=|g(\bm{X}(q),\bm{Y}(q))|\leq |\bm{X}(q)|_g|\bm{Y}(q)|_g=|\bm{Y}(q)|_g=d_0,
\]
which is a contradiction. Hence $L\subseteq N_0$ and consequently picking some $p\in L$ we obtain from the above identity $c_0=|\bm{Y}(p)|_g=d_0$. Therefore, we find for any $q\in N_0$
\[
|\bm{X}(q)|_g|\bm{Y}(q)|_g=d_0=c_0=g(\bm{Y}(q),\bm{X}(q)),
\]
so that by means of the Cauchy-Schwarz inequality $\bm{X}(q)$ and $\bm{Y}(q)$ must be linearly dependent, i.e. $N_0\subseteq L$. Overall $L=N_0$ and as pointed out before the proof is complete. $\square$
\begin{rem}[Topology of optimal domains]
	\label{S5R3}
To conclude the discussion of the topology of potential optimal domains, let $(\mathcal{R},g)$ be an ambient space (not necessarily with the div-free extension property). Then \cref{S2L6} still applies (only the constant may change values on distinct boundary components, which must necessarily all be non-zero if $\bm{X}$ is not the zero vector field by means of \cite[Lemma 2.1]{G21a}), so that all boundary components of any given optimal domain $\bar{M}$ must be tori. Then \cite[Proposition 18.6.2]{Die08} implies
\[
0=\chi\left(\partial\bar{M}\right)=2\chi\left(\bar{M}\right)
\]
and consequently, since the boundary is non-empty
\begin{align}
	\label{S5E3}
	\dim(H^1_{dR}(\bar{M}))=1+\dim(H^2_{dR}(\bar{M})).
\end{align}
If now in addition $\mathcal{R}=\mathbb{R}^3$ or $\mathcal{R}=B_1(0)\subset \mathbb{R}^3$ (but the metric $g$ need not be related to the Euclidean metric in any way), then \cite{CDG02} implies that $\dim(H^2_{dR}(\bar{M}))=\#\partial\bar{M}-1$ and therefore the de Rham cohomology groups of any given optimal domain are uniquely determined by the number of connected components of $\partial\bar{M}$, each of which is a torus.
\end{rem}
\section{Proof of \cref{MC8}}
Let us now come to the proof of \cref{MC8}.
\newline
\newline
\underline{Proof of \cref{MC8}:} Suppose for a contradiction that $\bar{M}$ is optimal. We use the notation preceding the statement of \cref{MC8}, recall also \cref{Fig1}. It then follows from \cref{MT2} that the eigenvalue $\lambda_+(\bar{M})$ must be simple. We can then argue identically as in the proof of \cref{S5P1} up until (\ref{S5E2}) to conclude the same equation and that $\bar{M}$ does not intersect the $z$-axis. Now by simplicity of the eigenvalue, the vector field $\bm{R}$, as it induces isometries and is tangent to the boundary, must commute with any fixed eigenfield $\bm{X}\in \mathcal{V}_n(\bar{M})$ corresponding to $\lambda_+(\bar{M})$. The result \cite[Lemma 3.5]{G21b} implies that we can find a vector potential $\bm{A}\in \mathcal{V}(\bar{M})$ of $\bm{X}$ which is normal to the boundary and commutes with $\bm{R}$ as well. Hence, we conclude from the Killing equations
\[
\text{grad}(g_E(\lambda_+(\bar{M})\bm{A},\bm{R}))=\lambda_+(\bar{M})[\bm{R},\bm{A}]+\lambda_+(\bar{M})\bm{R}\times \text{curl}(\bm{A})=\lambda_+(\bar{M})\bm{R}\times \bm{X}.
\]
Replacing $\lambda_+(\bar{M})\bm{A}$ by $\bm{X}$ and since $\bm{X}$ commutes with $\bm{R}$ and is an eigenfield of curl, we obtain in the same way
\[
\text{grad}(g_E(\bm{X},\bm{R}))=\lambda_+(\bar{M})\bm{R}\times \bm{X},
\]
so that we find $-\lambda_+(\bar{M})g_E(\bm{A},\bm{R})+g_E(\bm{X},\bm{R})=c_0$ for some $c_0\in \mathbb{R}$ on $\bar{M}$. But since $\bm{A}$ is normal to the boundary, while $\bm{R}$ is tangent to it, we obtain
\begin{align}
	\label{S6E4}
	g_E(\bm{X},\bm{R})=c_0\text{ on }\partial\bar{M}.
\end{align}
The important observation here is that the constant $c_0$ in (\ref{S6E4}) is the same on each boundary component. Now let $\partial \bar{M}_i$, $i=0,\dots,\#\partial\bar{M}-1=:N$, denote the boundary components of $\partial\bar{M}$ induced by the curves $\gamma_i$ respectively. Since $\partial\bar{M}_0$ is the boundary component containing the closest points to the $z$-axis it follows from the Cauchy-Schwarz inequality and (\ref{S6E4}) that $\bm{X}$ and $\bm{R}$ are, in a given cross section, linearly dependent, if at all, precisely on the set $N_0$, see (\ref{ME7}). Let us now as in the proof of \cref{S5P1} define $\bm{\Gamma}:=\frac{\bm{R}}{|\bm{R}|^2_{2}}\in \mathcal{H}_N(\bar{M})$, then since $\bm{X}$ and $\bm{\Gamma}$ are tangent to the boundary, we must have $\bm{X}\times \bm{\Gamma}=\pm |\bm{X}\times \bm{\Gamma}|_{2}\mathcal{N}$ at each boundary point (the sign may differ at different points), where $\mathcal{N}$ as usual is the outward unit normal. Therefore we have at every boundary point $|g_E(\bm{X}\times \bm{\Gamma},\mathcal{N})|=|\bm{X}\times \bm{\Gamma}|_{2}$. Further, we note that since $\bm{X}$ and $\bm{\Gamma}$ are, in a given cross section, linearly dependent at most on $N_0$, that after possibly replacing $\bm{X}$ by $-\bm{X}$ we find $\bm{X}\times \bm{\Gamma}=+|\bm{X}\times \bm{\Gamma}|_2\mathcal{N}$ on $\gamma_0\setminus L_-$. Hence we conclude from (\ref{S5E2}) that
\[
0=\int_{\partial\bar{M}}g_E(\bm{X}\times \bm{\Gamma},\mathcal{N})\omega_{g_{\partial\bar{M}}} 
\]
\[
\Leftrightarrow \int_{\left(\gamma_0\setminus L_-\right)\times S^1}|\bm{X}\times \bm{\Gamma}|_{2}\omega_{g_{\partial\bar{M}}}=\int_{L_-\times S^1}\pm |\bm{X}\times \bm{\Gamma}|_{2}\omega_{g_{\partial\bar{M}}}+\sum_{i=1}^N\int_{\partial\bar{M}_i} \pm|\bm{X}\times \bm{\Gamma}|_{2}\omega_{g_{\partial\bar{M}}},
\]
where the $\pm$ may be different for distinct points. Thus the triangle inequality implies
\begin{align}
	\label{S6E5}
	\int_{\left(\gamma_0\setminus L_-\right)\times S^1}|\bm{X}\times \bm{\Gamma}|_{2}\omega_{g_{\partial\bar{M}}}\leq \int_{L_-\times S^1}|\bm{X}\times \bm{\Gamma}|_{2}\omega_{g_{\partial\bar{M}}}+\sum_{i=1}^N\int_{\partial\bar{M}_i} |\bm{X}\times \bm{\Gamma}|_{2}\omega_{g_{\partial\bar{M}}}.
\end{align}
After an appropriate scaling we can achieve $|\bm{X}|_{2}=1$ on $\partial\bar{M}$ and then compute
\[
|\bm{X}\times \bm{\Gamma}|^2_{2}=|\bm{\Gamma}|^2_{2}-g^2_E(\bm{X},\bm{\Gamma})=\frac{1}{|\bm{R}|^2_{2}}-\frac{c^2_0}{|\bm{R}|^4_{2}}\text{ on }\partial\bar{M},
\]
where we used the definition of $\bm{\Gamma}$ and (\ref{S6E4}). Inserting this in (\ref{S6E5}) yields
\begin{align}
	\label{S6E6}
	\int_{\left(\gamma_0\setminus L_-\right)\times S^1}\frac{\sqrt{1-\frac{c^2_0}{|\bm{R}|^2_{2}}}}{|\bm{R}|_{2}}\omega_{g_{\partial\bar{M}}}\leq \int_{L_-\times S^1}\frac{\sqrt{1-\frac{c^2_0}{|\bm{R}|^2_{2}}}}{|\bm{R}|_{2}}\omega_{g_{\partial\bar{M}}}+\sum_{i=1}^N\int_{\partial\bar{M}_i} \frac{\sqrt{1-\frac{c^2_0}{|\bm{R}|^2_{2}}}}{|\bm{R}|_{2}}\omega_{g_{\partial\bar{M}}}.
\end{align}
Since we are dealing with surfaces of revolutions and the boundary components $\partial\bar{M}_i$ are induced by the closed curves $\gamma_i$ (not intersecting the $z$-axis), we see that for any function $f\in C^{\infty}(\partial\bar{M})$ with $\mathcal{L}_{\bm{R}}(f)=0$ ($\mathcal{L}$ denoting the Lie-derivative), we have
\[
\int_{\partial\bar{M}_i}f\omega_{g_{\partial\bar{M}}}=2\pi \int_{0}^{\mathcal{L}^1(\gamma_i)}\gamma^1_i(t)f(\gamma_i(t))dt=2\pi \int_0^{\mathcal{L}^1(\gamma_i)}|\bm{R}(\gamma_i(t))|_{2}f(\gamma_i(t))dt\text{ for all }1\leq i\leq N,
\]
where we without loss of generality assume that our curves are parametrised by arc-length, we consider a section in the $x$-$z$-plane which lies entirely in the first quadrant (which can always be arranged by translating the domain along the $z$-axis) and where $\gamma^1_i$ denotes the projection on the $x$-component of $\gamma_i$. A similar reasoning applies to the integrals with domain $(\gamma_0\setminus L_-)\times S^1$ and $L_-\times S^1$. Inserting this in (\ref{S6E6}) yields
\begin{align}
	\label{S6E7}
	\int_0^{\mathcal{L}^1(\gamma_0\setminus L_-)}\sqrt{1-\frac{c^2_0}{|\bm{R}(\gamma_0(t))|^2_{2}}}dt\leq \int_0^{\mathcal{L}^1( L_-)}\sqrt{1-\frac{c^2_0}{|\bm{R}(\gamma_0(t))|^2_{2}}}dt+\sum_{i=1}^N\int_0^{\mathcal{L}^1(\gamma_i)}\sqrt{1-\frac{c^2_0}{|\bm{R}(\gamma_i(t))|^2_{2}}}dt,
\end{align}
where we from now on always reparametrise parts of $\gamma_0$ by appropriate translations in 'time' so that the integrals all start at '$t=0$'. We show now that (\ref{S6E7}) contradicts our assumptions. First suppose that $c_0=0$, then our assumptions in combination with (\ref{S6E7}) imply
\begin{align}
	\label{S6E8}
\mathcal{L}^1(L_-)+\sum_{i=1}^N\mathcal{L}^1(\gamma_i)\leq \mathcal{L}^1(L_+)\leq \mathcal{L}^1(\gamma_0\setminus L_-)\leq \mathcal{L}^1(L_-)+\sum_{i=1}^N\mathcal{L}^1(\gamma_i),
\end{align}
so that $\mathcal{L}^1(L_+)=\mathcal{L}^1(\gamma_0\setminus L_-)$, which by definition of $L_+$ and $L_-$ then implies that $d_+=d_-$ as $L_+$ in this case must connect the points $\vec{x}_-$ and $\vec{x}_+$ which have minimal distance to the $z$-axis. This in particular implies that $N=0$ and by definition of $d_+$ the curve $L_-$ must be a line segment parallel to the $z$-axis (possibly just a point) which connects $\vec{x}_-$ and $\vec{x}_+$. Since $L_+$ connects the same points and the curve $\gamma_0$ is an embedded circle, we must have $\mathcal{L}^1(L_+)>\mathcal{L}^1(L_-)$. This contradicts (\ref{S6E8}).
\newline
So from now on let $c_0\neq 0$. We note that since on $L_+$ we have $|\bm{R}|_{2}\geq d_+$ and on $L_-$ and the $\gamma_i$ we have $|\bm{R}|_{2}\leq d_+$, we find
\[
\sqrt{1-\frac{c^2_0}{|\bm{R}(\gamma_0(t))|^2_{2}}}\geq \sqrt{1-\frac{c^2_0}{d^2_+}}\text{ on }L_+,
\]
\[
\sqrt{1-\frac{c^2_0}{|\bm{R}(\gamma_0(t))|^2_{2}}}\leq \sqrt{1-\frac{c^2_0}{d^2_+}}\text{ on }L_-\text{ and }\gamma_i\text{, }1\leq i\leq N.
\]
We observe that for at least one point on $L_+$ the above inequality must be strict, since otherwise $L_+$ consists of line segments (possibly points) parallel to the $z$-axis at distance $d_+$ from the axis. Hence in that case $L_+$ must necessarily intersect one of the $\gamma_i$ or $L_-$ at some point at which the distance $d_+$ is realised, which cannot happen since we are dealing with a disjoint collection of embedded circles. Then by continuity and our assumption we conclude
\[
\int_0^{\mathcal{L}^1(L_+)}\sqrt{1-\frac{c^2_0}{|\bm{R}(\gamma_0(t))|^2_{2}}}dt> \mathcal{L}^1(L_+)\sqrt{1-\frac{c^2_0}{d^2_+}}\geq \mathcal{L}^1(L_-)\sqrt{1-\frac{c^2_0}{d^2_+}}+\sum_{i=1}^N\mathcal{L}^1(\gamma_i)\sqrt{1-\frac{c^2_0}{d^2_+}}
\]
\[
\geq \int_0^{\mathcal{L}^1(L_-)}\sqrt{1-\frac{c^2_0}{|\bm{R}(\gamma_0(t))|^2_{2}}}dt+\sum_{i=1}^N\int_0^{\mathcal{L}^1(\gamma_i)}\sqrt{1-\frac{c^2_0}{|\bm{R}(\gamma_i(t))|^2_{2}}}dt.
\]
Estimating the first integral above by the left integral in (\ref{S6E7}) we obtain a contradiction. Hence in any case (\ref{S6E7}) contradicts our assumptions, which concludes the proof. $\square$
\section*{Acknowledgements}
This work has received funding from the European Research Council (ERC) under the European Union’s Horizon 2020 research and innovation programme through the grant agreement 862342.
\newline
Further, this work has been partially supported by the grant CEX2019-000904-S funded by MCIN\slash AEI/10.13039/501100011033. Parts of this work have also been supported by the Deutsche Forschungsgemeinschaft (DFG, German Research Foundation) – Projektnummer 320021702/GRK2326 –  Energy, Entropy, and Dissipative Dynamics (EDDy).
\appendix
\section{Killing-Beltrami fields and geometry}
In this section we provide some applications of a second variation inequality related to our isoperimetric problem. Since the focus of the present work lies on the model spaces of constant sectional curvature, we will formulate two results, which establish a connection between the existence of Killing-Beltrami fields, i.e. vector fields which satisfy the Killing equations and are eigenfields of the curl operator, and the sectional curvature.
\begin{thm}[Killing-Beltrami fields and sectional curvature I]
	\label{MT9}
	Let $(\mathcal{R},g)$ be an oriented, connected, smooth Riemannian $3$-manifold with (possibly empty) boundary and of constant sectional curvature $c\in \mathbb{R}$. Further let $\bm{Y}\in \mathcal{V}(\mathcal{R})\setminus \{0\}$ satisfy the Killing equations as well as $\text{curl}(\bm{Y})=\lambda \bm{Y}$ for some $\lambda\in \mathbb{R}$. Then
	\begin{align}
		\label{ME9}
		c=\frac{\lambda^2}{4}.
	\end{align}
	In particular, the only Killing field on a Riemannian manifold of constant negative sectional curvature satisfying $\text{curl}(\bm{Y})=\lambda \bm{Y}$ for some constant $\lambda$ is the zero vector field.
\end{thm}
Note that we do not demand in \cref{MT9} that $\bm{Y}\in \mathcal{V}_n(\bar{M})$ and we specifically allow $\lambda$ to be zero (in contrast to our definition of a (strong) Beltrami field). Further, we do not require the underlying metric to be complete.
\newline
As an example of \cref{MT9} we may consider the coordinate vector fields $e_i$, $i=1,2,3$, on Euclidean space or the standard Hopf vector fields on $S^3$ which are Killing fields, corresponding to the curl eigenvalue $\lambda=2$.
\newline
\newline
The following presents a situation in which the existence of Killing-Beltrami fields implies that the sectional curvature must be constant.
\begin{thm}[Killing-Beltrami fields and sectional curvature II]
	\label{MT10}
	Let $(\mathcal{R},g)$ be an oriented, connected, smooth Riemannian $3$-manifold with (possibly empty) boundary. Further, suppose that there exist two vector fields $\bm{Y}_1,\bm{Y}_2\in \mathcal{V}(\mathcal{R})$, satisfying the Killing equations, and a constant $\lambda\in \mathbb{R}$ with $\text{curl}(\bm{Y}_i)=\lambda \bm{Y}_i$ for $i=1,2$. If $\bm{Y}_1(p)\times \bm{Y}_2(p)\neq 0$ at some point $p\in \mathcal{R}$, then $(\mathcal{R},g)$ is of constant sectional curvature $c=\frac{\lambda^2}{4}$.
\end{thm}
Simple examples of \cref{MT10} are the Euclidean $3$-space with the standard coordinate vector fields $e_i$, $i=1,2,3$ and the flat $3$-torus $T^3$ (viewed as a cube in $\mathbb{R}^3$ with opposite faces identified) with the same vector fields (which descend to well-defined vector fields on $T^3$) or once more $S^3$ with its Hopf fields.
\newline
\newline
Both theorems will turn out to be a consequence of the following result
\begin{lem}
	\label{S7L1}
	Let $(\mathcal{R},g)$ be an oriented, connected, smooth Riemannian $3$-manifold with (possibly empty) boundary and let $\bm{Y}\in \mathcal{V}(\mathcal{R})$ be a Killing field, i.e. it satisfies the Killing equations. Then for every $\bm{X}\in \mathcal{V}(\mathcal{R})$ we have the following identity
	\[
	\text{sec}(\bm{X},\bm{Y})|\bm{X}\wedge \bm{Y}|^2_g=|\nabla_{\bm{X}}\bm{Y}|^2_g+g(\bm{X},\nabla_{\bm{X}}\nabla_{\bm{Y}}\bm{Y}),
	\]
	where $\text{sec}$ denotes the sectional curvature (and is set to zero by convention whenever the two vector fields are linearly dependent), $\nabla$ denotes the Levi-Civita connection and $|\bm{X}\wedge \bm{Y}|^2_g=|\bm{X}|^2_g|\bm{Y}|^2_g-g^2(\bm{X},\bm{Y})$.
\end{lem}
Before we come to the proof of \cref{S7L1} let us shortly discuss its main idea and its relation to the isoperimetric problem we studied so far. Namely, in order to derive \cref{S2L6} one considers some fixed optimal domain $\bar{M}$ and eigenfield $\bm{X}\in \mathcal{V}_n(\bar{M})$ of $\lambda_+(\bar{M})$. Then for any fixed divergence-free vector field $\bm{W}\in \mathcal{V}(\bar{M})$ we let $\psi_t$ denote the induced flow of some divergence-free extension $\widetilde{\bm{W}}$ of $\bm{W}$ to some open neighbourhood $U$ of $\bar{M}$. One then defines the following function for suitable $0<\epsilon\ll 1$
\begin{align}
	\label{S7E1}
	f:(-\epsilon,\epsilon)\rightarrow \mathbb{R},\text{ }t\mapsto \int_{\psi_t(\bar{M})}|(\psi_t)_{*}\bm{X}|^2_g\omega_g,
\end{align}
where $(\psi_t)_{*}\bm{X}$ denotes the pushforward vector field. It turns out that this function must have a global minimum at $t=0$, so that we obtain the necessary conditions $\frac{d}{dt}|_{t=0}f(t)=0$ and $\frac{d^2}{dt^2}|_{t=0}f(t)\geq 0$. The first equation then leads to \cref{S2L6}, see \cite[Theorem D]{CDGT002} and \cite[Theorem 2.2.3 and Proposition 2.2.4]{G20Diss}. Further, one can compute the second variation to obtain an integral inequality in the setting of the isoperimetric problem, which was done in \cite[Theorem 2.2.5]{G20Diss}. Now if the vector field $\bm{W}$ is not only divergence-free but in fact extends to a Killing field on $\mathcal{R}$, then its flow acts by isometries and hence the $L^2$-energy is preserved under its action, i.e. the function $f$ defined in (\ref{S7E1}) is independent of $t$. This remains true even if $\bar{M}$ is not an optimal domain and if $\bm{X}$ is any smooth vector field on $\bar{M}$. Hence the second variation leads to an integral identity, rather than an inequality in this case, which can be averaged.
\newline
\newline
\underline{Proof of \cref{S7L1}:} We may restrict to the interior of $\mathcal{R}$ by means of a density argument. So let $p\in \text{int}(\mathcal{R})$ be any fixed point and let $0<r\ll 1$ be fixed, such that the closed geodesic ball $B_r(p)$ is contained in $\text{int}(\mathcal{R})$. Now if $\bm{Y}$ is our Killing field it induces a flow of isometries, so that we may consider the function $f$ as defined in (\ref{S7E1}) with $\bar{M}=B_r(p)$. As explained in the discussion following (\ref{S7E1}) the function $f$ is independent of $t$, so that we must have $\frac{d^2}{dt^2}|_{t=0}f(t)=0$. Doing the computations, see the proof of \cite[Theorem 2.2.5, Equation (2.10.14)]{G20Diss} for the details, yields the following
\[
\int_{B_r(p)}\text{sec}(\bm{X},\bm{Y})|\bm{X}\wedge \bm{Y}|^2_g\omega_g=\int_{B_r(p)}|\nabla_{\bm{X}}\bm{Y}|^2_g\omega_g+\int_{B_r(p)}g(\bm{X},\nabla_{\bm{X}}\nabla_{\bm{Y}}\bm{Y})\omega_g.
\]
Now averaging and taking the limit $r\searrow 0$, continuity gives the desired result. $\square$
\newline
\newline
\underline{Proof of \cref{MT9}:} Let $\bm{Y}$ be a non-trivial Killing field on $\mathcal{R}$ satisfying $\text{curl}(\bm{Y})=\lambda \bm{Y}$. We observe that the Killing equations imply
\begin{align}
\label{S7E2}
\nabla_{\bm{X}}{\bm{Y}}=\frac{\bm{X}\times \text{curl}(\bm{Y})}{2}=\frac{\lambda}{2}\bm{X}\times \bm{Y}\text{ for every }\bm{X}\in \mathcal{V}(\mathcal{R}),
\end{align}
the latter by means of the Beltrami property. Since $\bm{Y}$ is non-trivial we can find some point $p$ at which $\bm{Y}(p)\neq 0$ (we may assume that $p$ is an interior point). Now fix any tangent vector $X\in T_p\mathcal{R}$ which is linearly independent from $\bm{Y}(p)$ and extend it to a smooth vector field $\bm{X}\in \mathcal{V}(\mathcal{R})$ with $\bm{X}(p)=X$. We note that (\ref{S7E2}) implies, by setting therein $\bm{X}=\bm{Y}$, that $\nabla_{\bm{Y}}\bm{Y}=0$, so that it follows from \cref{S7L1} and (\ref{S7E2}) that
\[
c|\bm{X}\wedge\bm{Y}|^2_g=\text{sec}(\bm{X},\bm{Y})|\bm{X}\wedge\bm{Y}|^2_g=\frac{\lambda^2}{4}|\bm{X}\times \bm{Y}|^2_g=\frac{\lambda^2}{4}|\bm{X}\wedge\bm{Y}|^2_g.
\]
We note that $|\bm{X}\wedge\bm{Y}|^2_g=0$ at some point if and only if $\bm{X}$ and $\bm{Y}$ are linearly dependent at that point, so that we find $c=\frac{\lambda^2}{4}$. $\square$
\newline
\newline
Before stating the proof of \cref{MT10} we will establish the following unique continuation result
\begin{lem}
\label{S7L2}
Let $(\mathcal{R},g)$ be an oriented, connected, smooth Riemannian $3$-manifold with (possibly empty) boundary. Further let $\bm{X}_i\in \mathcal{V}(\mathcal{R})$, $i=1,2$, be two smooth divergence-free vector fields such that $\text{curl}(\bm{X}_i)=\lambda\bm{X}_i$, $i=1,2$, for some constant $\lambda\in \mathbb{R}$. If the set $\{\bm{X}_1\times \bm{X}_2=0\}$ has an interior point, then $\bm{X}_1=\kappa\bm{X}_2$ for some constant $\kappa\in \mathbb{R}$.
\end{lem}
\underline{Proof of \cref{S7L2}:} Let $W\subseteq \{\bm{X}_1\times \bm{X}_2=0\}$ be a non-empty open subset. It follows from \cite{AKS62} that the sets $W_i:=\{\bm{X}_i\neq 0\}$, $i=1,2$, are open and dense in $\mathcal{R}$. After replacing $W$ by $W\cap W_1\cap W_2$ and shrinking it further if necessary we obtain a non-empty, connected open set $W$ on which the $\bm{X}_i$ are non-vanishing and linearly dependent. We conclude that there exists a smooth function $f\in C^{\infty}(W)$ with $\bm{X}_1=f\bm{X}_2$. Now by the eigenfield property we find
\[
\lambda \bm{X}_1=\text{curl}(\bm{X}_1)=\text{curl}(f\bm{X}_2)=\text{grad}(f)\times \bm{X}_2+f\lambda \bm{X}_2=\text{grad}(f)\times \bm{X}_2+\lambda \bm{X}_1
\]
and consequently $\text{grad}(f)\times \bm{X}_2=0$. Since $\bm{X}_2$ is nowhere vanishing on $W$ we find a smooth function $h\in C^{\infty}(W)$ with $\text{grad}(f)=h\bm{X}_2$. Now since $\bm{X}_1$ and $\bm{X}_2$ are divergence-free we conclude
\[
0=\text{div}(\bm{X}_1)=\text{div}(f\bm{X}_2)=hg(\bm{X}_2,\bm{X}_2),
\]
where we used that $\text{grad}(f)=h\bm{X}_2$. Since $\bm{X}_2$ is nowhere vanishing, we must have $h=0$ on $W$ and hence $\text{grad}(f)=0$ on $W$, so that $f=\kappa\in \mathbb{R}$ is some constant with $\bm{X}_1=\kappa \bm{X}_2$ on $W$. We conclude that the vector field $\bm{X}_1-\kappa\bm{X}_2$ is divergence-free and a curl eigenfield corresponding to the eigenvalue $\lambda$, which vanishes identically on some non-empty open subset. Then \cite{AKS62} implies that this vector field must be identically zero on all of $\mathcal{R}$. $\square$
\newline
\newline
\underline{Proof of \cref{MT10}:} Let $\bm{Y}_1$ and $\bm{Y}_2$ be two smooth Killing fields on $\mathcal{R}$, which are linearly independent in at least one point and which satisfy $\text{curl}(\bm{Y}_i)=\lambda \bm{Y}_i$, $i=1,2$, for the same constant $\lambda$. Since Killing fields are divergence-free, we conclude from \cref{S7L2} that the set $\mathcal{U}:=\{\bm{Y}_1\times \bm{Y}_2\neq 0\}$ is an open and dense subset of $\mathcal{R}$. It then follows that $\left(\bm{Y}_1,\bm{Y}_2,\bm{Y}_1\times \bm{Y}_2\right)$ gives rise to a smooth frame of the tangent bundle of $\mathcal{R}$ on $\mathcal{U}$. We conclude from an identical argument as in the proof of \cref{MT9} that
\[
\text{sec}(\bm{Y}_1,\bm{Y}_2)=\frac{\lambda^2}{4}\text{ on }\mathcal{U}.
\]
Now if we fix any $p\in \mathcal{U}$ and any two linearly independent tangent vectors $X,Z\in T_p\mathcal{R}$, then we can express them as $X=\sum_{i=1}^2X^i\bm{Y}_i+X^3\bm{Y}_1\times \bm{Y}_2$ and $Z=\sum_{i=1}^2Z^i\bm{Y}_i+Z^3\bm{Y}_1\times \bm{Y}_2$. If $X^3=0=Z^3$, then since the sectional curvature solely depends on the space spanned by the tangent vectors involved, we find $\text{sec}(X,Z)=\text{sec}(\bm{Y}_1,\bm{Y}_2)=\frac{\lambda^2}{4}$. We are left with considering the case $X^3\neq 0$. We may assume $X^3=1$ and $Z^3=0$, since the sectional curvature only depends on the space spanned by $X$ and $Z$. We therefore conclude that $\text{sec}(X,Z)=\text{sec}(X,\bm{Y})$, where $\bm{Y}=\sum_{i=1}^2Z^i\bm{Y}_i$. But $\bm{Y}$ is now itself a Killing field and an eigenfield of curl corresponding to the eigenvalue $\lambda$, so that once more the reasoning of the proof of \cref{MT9} applies and we conclude $\text{sec}(X,Z)=\text{sec}(X,\bm{Y})=\frac{\lambda^2}{4}$. Hence for any choice of $p\in \mathcal{U}$ and any linearly independent vectors $X,Z\in T_p\mathcal{R}$ we have $\text{sec}(X,Z)=\frac{\lambda^2}{4}$. Now lastly if we let $p\in \mathcal{R}\setminus \mathcal{U}$ and $X,Z\in T_p\mathcal{R}$ be any linearly independent vectors, then we may extend them to smooth, linearly independent vector fields $\bm{X}$, $\bm{Z}$ to some open neighbourhood $W$ of $p$ with $\bm{X}(p)=X$ and $\bm{Z}(p)=Z$. We can then find a sequence $(p_k)_k\subset \mathcal{U}$ converging to $p$ and compute the sectional curvature as follows \cite[Proposition 8.29]{L18}, where $\text{Rm}$ denotes the Riemann curvature tensor
\[
\text{sec}(X,Z)=\text{sec}(\bm{X}(p),\bm{Z}(p))=\frac{\text{Rm}(\bm{X}(p),\bm{Z}(p),\bm{Z}(p),\bm{X}(p))}{|\bm{X}(p)\wedge \bm{Z}(p)|^2_g}
\]
\[
=\lim_{k\rightarrow\infty}\frac{\text{Rm}(\bm{X}(p_k),\bm{Z}(p_k),\bm{Z}(p_k),\bm{X}(p_k))}{|\bm{X}(p_k)\wedge \bm{Z}(p_k)|^2_g}=\lim_{k\rightarrow\infty}\text{sec}(\bm{X}(p_k),\bm{Z}(p_k))=\frac{\lambda^2}{4}.
\]
This concludes the proof. $\square$
\section{Isometric embeddings and flat surfaces}
\begin{lem}
\label{AL1}
Let $(\mathcal{R},g)\in \{\mathbb{R}^3,\mathcal{H}^3,S^3_+\}$ be our ambient space. If $S\subset \mathcal{R}$ is any (oriented) closed hypersurface, then there exists some $p\in S$ such that $\kappa(p)\neq 0$, where $\kappa:S\rightarrow \mathbb{R}$ denotes the induced Gaussian curvature.
\end{lem}
\underline{Proof of \cref{AL1}:} We consider as usual the Poincar\'{e} ball model for $\mathcal{H}^3$ and identify $S^3_+$ with $B_1(0)$ equipped with the metric $g(x):=\frac{4}{(1+|x|^2_2)^2}g_E(x)$, where $g_E$ is the Euclidean metric on the open unit ball $B_1(0)$. First we note that by compactness of $S$ there exists some point $p\in S$ such that $|\cdot|^2_2:S\rightarrow \mathbb{R}$ admits a global maximum at $p$. Now let $E_1\in T_pS$ be a normalised eigenvector of the shape operator of $S$ at $p$. We can then consider a unit speed geodesic of $S$, $\gamma:(-\epsilon,\epsilon)\rightarrow S\subset \mathcal{R}$ with $\gamma(0)=p$ and $\dot{\gamma}(0)=E_1$. Then we find
\[
\kappa_1(p)\mathcal{N}(p)=\nabla_{\dot{\gamma}}{\dot{\gamma}}(0),
\]
where $\kappa_1(p)$ denotes the eigenvalue corresponding to $E_1$, $\nabla_{\dot{\gamma}}{\dot{\gamma}}$ denotes the acceleration of $\gamma$ viewed as a curve into $\mathcal{R}$ and $\mathcal{N}$ is the unit normal with respect to which the shape operator is considered, see \cite[Proposition 8.10]{L18}. Thus, if we can show that the corresponding curvature of $\gamma$ at $t=0$ is non-zero, then $\kappa_1(p)\neq 0$. To see that this is the case we define $f(t):=\frac{|\gamma(t)|^2_2}{2}$ and note that $f$ has a global maximum at $t=0$ by choice of $p$. Hence
\begin{equation}
	\label{AE1}
0=\frac{d}{dt}|_{t=0}f(t)=\gamma(0)\cdot \dot{\gamma}(0)\text{ and }0\geq \frac{d^2}{dt^2}|_{t=0}f(t)=|\dot{\gamma}(0)|_2^2+\gamma(0)\cdot \ddot{\gamma}(0),
\end{equation}
where $\cdot$ denotes the standard Euclidean inner product. If we consider a conformal metric $g=\exp(2\phi)g_E$ for some smooth function $\phi:B_1(0)\rightarrow \mathbb{R}$ on the open unit ball, then the Levi-Civita connection changes according to the following well-known formula
\[
\nabla^g_{\dot{\gamma}}\dot{\gamma}(0)=\ddot{\gamma}(0)+2(\text{grad}(\phi)(\gamma(0))\cdot \dot{\gamma}(0))\dot{\gamma}(0)-|\dot{\gamma}(0)|^2_2\text{grad}(\phi)(\gamma(0)),
\]
where all quantities on the right hand side are computed with respect to the Euclidean metric. Multiply this identity by $\gamma(0)$ and use the first identity from (\ref{AE1}) to find
\[
0\geq |\dot{\gamma}(0)|^2_2\left(1+\gamma(0)\cdot \text{grad}(\phi)(\gamma(0))\right)+\nabla^g_{\dot{\gamma}}\dot{\gamma}(0)\cdot \gamma(0)
\]
\[
\geq |\dot{\gamma}(0)|^2_2\left(1+\gamma(0)\cdot \text{grad}(\phi)(\gamma(0))\right)-|\nabla^g_{\dot{\gamma}}\dot{\gamma}(0)|_2|\gamma(0)|_2
\]
and so $|\nabla^g_{\dot{\gamma}}\dot{\gamma}(0)|_2\geq \frac{|\dot{\gamma}(0)|^2_2\left(1+\gamma(0)\cdot \text{grad}(\phi)(\gamma(0))\right)}{|\gamma(0)|_2}$. In particular $\nabla^g_{\dot{\gamma}}\dot{\gamma}(0)\neq 0$ as soon as $1+\gamma(0)\cdot \text{grad}(\phi)(\gamma(0))$ $>0$, which can be confirmed by direct calculations for $\mathbb{R}^3$, $\mathcal{H}^3$ and $S^3_+$, keeping in mind that $0<|\gamma(0)|_2<1$ for the non-Euclidean model spaces. Hence $\kappa_1(p)\neq 0$. In the same fashion we find for the second principle curvature $\kappa_2(p)\neq 0$ and overall $\kappa(p)=\kappa_1(p)\kappa_2(p)\neq 0$. $\square$
\bibliographystyle{plain}
\bibliography{mybibfileNOHYPERLINK}
\footnotesize
\end{document}